\documentclass[11pt]{article}
\input{epsf.sty}
\usepackage{hyperref}
\usepackage{amsmath, amssymb}
\newtheorem {thm}{Theorem}[section]
\newtheorem {prop}[thm]{Proposition} 
\newtheorem {fact}[thm]{Fact}
\newtheorem {lem}[thm]{Lemma}
\newtheorem {cor}[thm]{Corollary}
\newtheorem {defn}[thm]{Definition}

\newenvironment{remark}[1][Remark:]{\begin{trivlist}
\item[\hskip \labelsep {\bfseries #1}]}{\end{trivlist}}
\def\Cox{\hfill \Box}

\def\N{{\Bbb N}}

\def\R{{\Bbb R}}
\def\LL{{\Bbb L}}
\def\P{{\Bbb P}}

\def\E{{\Bbb E}}

\def\e{{\varepsilon}}

\def\ba{{\backslash}}

\def\D{\Delta}
\def\a{\alpha}

\def\ba{\setminus}
\def\b{\beta}
\def\MM{{\cal M}}
\def\d{\delta}

\def\e{\varepsilon}

\def\phi{\varphi}

\def\g{\gamma}

\def\td{\tilde}

\def\l{\lambda}

\def\r{\rho}
\def\LL{\mathcal{L}}
\def\s{\sigma}

\def\D{\Delta}

\def\O{{\Omega}}

\def\P{{\Phi}}

\def\T{\T}

\def\PP{{\cal P}}

\begin{document}
\title{Continuous Spin Mean-Field models: \\ 
Limiting kernels and Gibbs Properties of local transforms  \\ 
}

\author{
Christof K\"ulske
\footnote{ University of Groningen, Institute  of Mathematics and
Computing Science, Postbus 407, 9700  AK Groningen, The
Netherlands,
\texttt{kuelske@math.rug.nl}, \texttt{
http://www.math.rug.nl/$\sim$kuelske/ }}\, and
Alex A. Opoku
\footnote{ University of Groningen, Department of Mathematics and
Computing Science, Postbus 407, 9700 AK Groningen, The
Netherlands, \texttt{A.opoku@math.rug.nl}} 
}

\maketitle

\begin{abstract} 
We extend the notion of Gibbsianness for mean-field systems 
to the set-up of general (possibly continuous) local state spaces. 
We investigate the Gibbs properties of systems arising 
from an initial mean-field Gibbs measure by application of given local transition kernels. 
This generalizes previous case-studies made for spins taking finitely many values to 
the first step in direction to a general theory, containing the following parts:   
(1) A formula 
for the limiting conditional probability distributions of the 
transformed system. It holds both in the Gibbs and non-Gibbs regime  
and invokes a minimization problem for a "constrained rate-function". 
(2) A criterion for Gibbsianness of the transformed 
system for initial Lipschitz-Hamiltonians involving concentration properties 
of the transition kernels.  
(3) A continuity estimate for the single-site conditional distributions 
of the transformed system.  While (2) and (3) have provable lattice-counterparts, 
the characterization of (1) is stronger in mean-field.  
As applications we show short-time Gibbsianness of rotator mean-field models on the 
$(q-1)$-dimensional sphere under diffusive time-evolution and the preservation of Gibbsianness 
under local coarse-graining of the initial local spin space.

 \end{abstract}



\vfill\eject
\section{Introduction}

The study of the (failure of the) Gibbs property 
is a source of interesting probability theory 
and is linked to the study of phase-transitions.  
Gibbs measures and generalized Gibbs measures are 
of interest not only on the lattice, but on more general structures. 
Examples of such structures are random graphs, or, in the simplest 
conceivable case, the complete graph, where the models 
are called mean-field models.

The Gibbs property of a given measure 
should be viewed as a continuity property 
of conditional probabilities as a function of the conditioning.   
When one tries to prove or disprove this property  
for a measure obtained by an application of a deterministic or stochastic 
transformation from a well-understood initial measure 
one is led to a constrained (or "quenched") 
problem, with "quenched impurities" that are induced by the conditioning. 
This introduces a "random" (or in better words constrained) system  
that we need to understand \cite{ACD}, and this ties the problem 
to disordered systems and statistical mechanics on random structures. 
 
It was through rigorous implementation of 
Renormalization Group transformations  that it was discovered that  images of Gibbs 
measures can be non-Gibbs \cite{ACD,GP78,GP79,DEZ}. 
After this discovery, there has been an interest in recent times in particular 
in the study of the loss and possible  recovery of the Gibbs property of an initial Gibbs measure
under a stochastic time-evolution. 
The study started in \cite{ACD1} where the authors focused on the  
evolution of an initial Gibbs measure of a lattice spin
Ising model under high-temperature spin-flip Glauber dynamics.
 The main phenomenon observed here was the loss
of the Gibbs property after a certain transition time when the
system was started from an initial low temperature state. 
 The measure  stays non-Gibbs forever when the initial external field is zero.
More complicated transitions between Gibbs and non-Gibbs
are possible at intermediate times when there is no spin-flip
symmetry. The case of site-wise independent diffusions of continuous spins on the lattice
starting from the Gibbs-measure of  a special double-well potential
was considered in \cite{KUL5}, exhibiting similarities and differences to the Ising case. 
In \cite{WIO} the authors  studied  models for continuous
compact spins, namely the planar rotor models on the circle subjected
to diffusive time-evolution. It is shown therein that starting with 
an initial low-temperature Gibbs measure, the time-evolved measure 
obtained for infinite- or high-temperature dynamics stays Gibbs for 
short times and for the corresponding  initial infinite- or high-
temperature Gibbs measure under infinite- or high-temperature 
dynamics, the time-evolved measure stays Gibbs forever. Their 
analysis uses the machinery of cluster expansions, as earlier developed in 
\cite{ROE}. Even before it was shown that the whole process of space-time 
histories can be viewed as a Gibbs measure \cite{DEU} which however does not 
 imply that fixed-time projections are Gibbs. 

Let us move from concrete examples to the elements of a general theory 
which have been proved so far. 
In \cite{KULOP} the preservation of the Gibbs property for 
compact (discrete and continuous)  spin models for general initial  interactions 
(having a finite "triple-norm")  
subjected  to general site-wise transformations is studied. 
The technique employed therein is Dobrushin uniqueness \cite{DOB, GOR}, which is quite 
robust and gives rise to explicit estimates. We obtained both quantitative estimates 
on the parameter regimes where Gibbsianness provably holds and, as the main new part,  
explicit continuity estimates for the conditional probabilities of the transformed system. 

As an example it is shown therein that starting with an initial Gibbs 
measure of a rotator spin model 
on the $(q-1)$-dimensional sphere ($q\geq2$) and performing site-wise independent diffusive 
time-evolutions, the Gibbs property is preserved in an explicitly computable time-interval starting 
from zero. Similar conclusions were 
drawn for  Gibbs measures for general initial interactions (with compact metric local spin-spaces) 
subjected to a local coarse-graining transformation. 
(Given a decomposition of the local state-space $S$ into countably many sets, the corresponding 
local coarse-graining is the map that associates to any point in $S$ the label of the corresponding 
set in the decomposition.)
Here the Gibbs property is preserved 
whenever the diameter of the largest set in the decomposition 
is small enough.  Roughy speaking,  
this result can be seen as stability of Gibbsianness under application of a ball of 
sufficiently fine local transformations of coarse-graining type.


In a related line of research, transforms of initial Gibbs measures for various mean-field 
models were investigated.  A variety of measures has been found 
 to be non-Gibbs \cite{KUL2,KUL1,HaggKu04} in the mean-field sense. Usually the analysis 
 of such systems shows parallels to what can be done on the lattice but goes much further. 
  We remark that in all the cases studied so far, mean-field spins that take finitely many  values 
had been considered, and a unifying treatment including discrete and continuous spins 
had been lacking.  For state of the art reviews on 
Gibbsianness and non-Gibbsianness we refer the 
reader to \cite{ACD, AFE08}.

Now, in this note we present a systematic investigation 
of the Gibbs property of mean-field measures 
subjected to local kernels. We are out to extend previous 
results on spins taking finitely many values to general possibly continuous  (but compact) spins.
More mathematical care is needed since we consider distributions of empirical 
measures taking values in an infinite-dimensional space. 
So, let us provide an informal roadmap of the present paper now, leaving 
the precise definitions and statements of the theorems to the main body.

What are the initial measures we are dealing with? 
We start in Section \ref{s2}  by defining a class of interactions $\Phi(\nu)$ as functions 
on empirical measures $\nu$ of the system. 
The corresponding mean-field Hamiltonian in a volume of size $N$ is 
$N \Phi(\nu)$.  The densities 
of the finite-volume Gibbs-measures w.r.t. an a-priori product measure $\a$ in volume $N$ 
are given in terms of the normalized exponential $\frac{1}{Z_N}\exp(- N \Phi(\nu))$. 

The first decision to be made is to find an appropriate notion of regularity of allowed interactions $\Phi$. 
It turns out here that the natural requirement 
(the suitable mean-field analogue of  the standard notion of absolute summability for interactions on the 
lattice) is that of 
continuous differentiability (in the space of measures on the single-site configuration space).  

Next, we define the notion of mean-field Gibbsianness of a model 
which is given  in terms of the sequence of its finite-volume measures, by looking 
at large-volume limits of single-site conditional distributions obtained from this. 
This procedure provides 
us with a kernel  
$\g_1(d\s_i | \nu)$ where $\nu$ is the empirical measure of a configuration in the conditioning. 
The model is called Gibbs if  every $\nu$ is a continuity point of $\g$ (for the weak topology). 
This is  a natural generalization from the 
discrete spin-examples where this notion had been introduced and investigated before.   
From this definition it is also apparent that the regularity requirement on the interaction 
made above is natural since it implies Gibbsianness of the initial system (see more on this below 
Definition \ref{consistency}). 

The situation is easier (and thus amenable to analysis) 
for mean-field models than for lattice models, since a configuration in the conditioning 
is replaced by a measure $\nu$ on the local spin space, and it is just 
one single-site kernel that captures the large-volume behavior. 

In Section \ref{s3}
we turn to the main focus of the paper, namely two-layer models, obtained 
by applying a stochastic kernel, independently over the sites, to the initial model.
A complete analysis of the Ising model in mean-field 
under stochastic site-wise independent time-evolution has been given in \cite{KUL2}, showing 
the emergence of non-Gibbsianness at sharp critical times 
and a phenomenon called symmetry-breaking in the set 
of bad configurations. More examples are found in \cite{KUL1,HaggKu04}.  
At first we develop the general theory which relates our desired object, 
the large-system limiting conditional 
distribution of the transformed system, to a variational problem. 
In this part no specific assumptions (other than continuous differentiability of the initial potential) 
will be made on 
the model. The results hold in regions of the parameter space of the interaction 
where both Gibbsianness and non-Gibbsianness can occur. In the non-Gibbsian regime however
we have to stay away from the specific critical values of the conditionings 
for which non-unique global minimizers occur. 
For the convenience of the reader we briefly review some background material on large deviations 
we will use for our analysis. 
Large deviation principles are interesting in themselves, 
but from the point of view of this paper, they will just be used 
as a tool to treat the limiting conditional probabilities. 
The main general result of this first general part is 
Theorem \ref{maintheorem}  which describes the infinite-volume second-layer 
conditional probabilities in terms of a solution 
of a variational problem (leading to a consistency equation) 
for the constrained first layer model (CFLM). 

In Section \ref{s4} we provide criteria for 
Gibbsianness of the transformed model.  
This part is based on the study of the constrained 
consistency equation obtained in the first part of the paper.  By Tychonovs theorem there 
exists at least one solution. By the contraction mapping theorem 
there is precisely one solution, provided the respective kernel is Lipschitz, 
uniformly in the conditioning, 
with a constant $L$ which can be derived explicitly, when $L<1$. 
Uniqueness of the solution implies mean-field Gibbsianness 
of the transformed model, by the first part.  
This is in nice analogy to the corresponding lattice results obtained 
in the paper \cite{KULOP} 
using techniques based on Dobrushin
 uniqueness. 
More can be said however about the transformed system, and can be put 
in perspective with corresponding lattice results. 

In \cite{KULOP} we were proving Gibbsianness but we did more than that. 
We provided explicit continuity estimates of the form 
$$\Vert \g'_i(d\eta_i | \eta_{i^c})-\g'_i(d\eta_i | \eta'_{i^c})\Vert \leq \sum_{j: j\neq i} 
Q_{i,j}d'(\eta_j,\eta'_j)$$
where $\g'_i(d\eta_i |  \eta_{i^c}) $ are single-site conditional probabilities
for the transformed system, $Q_{ij}$ is the so-called Goodness-matrix, 
and $d'$ is the so-called posterior metric. The posterior metric is the variational 
distance between constrained single-site measures 
$d'(\eta_i,\eta'_i)=\Vert K(d\s_i|\eta_i)- K(d\s_i|\eta'_i)\Vert $ where 
$K$ is the joint single-site a priori measure (obtained in terms of $\a$ and the transformation 
kernel).   

In the present mean-field setup we prove as the main result of the second part of the 
paper an estimate of the form 
$$\Vert \g'_i(d\eta_i | \nu)-\g'_i(d\eta_i | \nu')\Vert \leq L_2 \Vert \nu-\nu'\Vert $$
with $L_2$ given in Theorem \ref{contestimates}. 
In the lattice estimate there is a matrix $Q$ appearing, describing the spatial decay 
of influence of a variation of the conditioning at site $j$ while in the mean-field 
estimate we are simply considering the variational distance of the empirical 
measure of the conditioning. 

$L_2$ will be finite  for an initial interaction that is arbitrarily large but Lipschitz when 
the constrained single-site measures have good concentration properties. 
This is the case e.g. at short-times for diffusive time-evolutions, 
or for sufficiently fine local coarse-grainings.   
When the initial interaction is small the transformation plays no role, 
and $L_2$ is finite always.

We conclude the paper with the discussion of 
stochastic time-evolutions and local coarse-grainings in Section \ref{s5}.

\bigskip 
\bigskip

\section{Generalities on Mean-Field Models}\label{s2}
\subsection{Set-up}
Let $(S,d)$ and $(S',d')$ be two given compact Polish spaces (compact separable  
metric spaces), each equipped 
with their corresponding Borel $\s$-algebras. We denote by $\PP(S)$, ${\MM}_+(S)$
and $\MM(S)$ ($\PP(S')$, ${\MM}_+(S')$ and $\MM(S')$) the spaces of probability 
measures, finite positive measures and finite signed measures on $S$ ($S'$) 
respectively. Let $\a$ and $\a'$ be two given reference Borel probability
measures (also called the {\em a priori measures}) on $S$ and $S'$ respectively . 
In the following we will refer to $S$ 
as the {\em initial (first-layer) single-site spin space} and $S'$ as also the 
{\em transformed (second-layer) single-site spin space}. We respectively write  $\O=S^\N$ and 
$\O'={S'}^\N$ as the configuration spaces for the initial (first-layer) and the 
transformed (second-layer) systems.
In the sequel we will write probability measures  for the transformed system with primes and those
for the joint system (comprising of the initial and transformed systems) with tildes.
The probability measures  for the initial system will always be written without primes and tildes. Again
we denote by $\s$, $\eta$ and $\xi$ the spin variable for the initial, the transformed and  the joint systems
respectively, (e.g $\xi=(\s_i,\eta_i)_{i\in\N}\in\td\O=(S\times S')^\N$).
We further  set $V_N=\{1,\cdots,N\}$ and write  $\s_{V_N}$ for 
points in the product space $S^N$. We will simply write $\s$ instead of $\s_\N$.
We now define the following  concept of mean-field interaction 
for the initial systems that we shall consider in this work.
\begin{defn}\label{mfinteraction}
We shall refer to a map $\P: {\MM}_+(S)\rightarrow \R$ as a proper mean-field
 interaction if it satisfies the following conditions:
\begin{enumerate}
\item it is weakly  continuous  and
\item it satisfies the uniform directional differentiability condition, 
meaning that, for each $\nu\in {\MM}_+(S)$ the derivative 
$\P^{(1)}(\nu,\mu)$ at $\nu$ in direction $\mu$ exists and we have 
\begin{equation}\label{differentiability}
\P(\nu+\mu)-\P(\nu)-\P^{(1)}(\nu,\mu)=r(\mu)
\end{equation}
with
 $\lim_{t\rightarrow 0^+}
\frac{r(t\mu)}{t}=0$ uniformly in  $\mu\in\MM(S)$ for which $\nu+t\mu\in\MM_+(S)$, 
for $t\in(0,1]$. 
\item $\P^{(1)}(\nu,\mu)$ is a continuous function of $\nu$ 
\end{enumerate}

\end{defn}




For each mean-field interaction $\P$ and each $N\in\N$ we define the finite-volume Hamiltonian 
$H_N$ (a real-valued function on the product space $S^N$) as
\begin{equation}\label{freebd}
H_N(\s_{V_N}):=N\P\big(L_N(\s_{V_N})\big),
\end{equation} 
where $L_{N}\left(
 \s_{V_N}\right)=\frac{1}{N}\sum_{i=1}^{N}\delta_{\sigma_{i}} $
is the empirical measure. Observe from the  permutation 
invariance of the empirical measures that $H_N$ is also permutation invariant.
With this notation we define the finite-volume  Gibbs measure $\mu_{\b,N}$
 for the finite-volume Hamiltonian $H_N$ and at inverse temperature $\beta$   as 
\begin{equation}\label{mean-field}
\mu_{\beta,N}(d\s_{V_N}):=\dfrac{e^{-\beta H_N\big(L_N(\s_{V_N})
 \big)}\a^{\otimes N}
(d\s_{V_N})}{\int_{E^N}e^{-\beta H_N\big(L_N(\bar\s_{V_N})
 \big)}\a^{\otimes N}
(d\bar\s_{V_N})},
\end{equation} 
where  we've used $\otimes$ to denote tensor product of measures.
 In the following, unless otherwise stated, the inverse temperature $\beta$  will be absorbed into the
interaction $\P$. In view of this, we will  write $\mu_N$ instead of 
$\mu_{\beta,N}$. It follows from the permutation invariance of the $H_N$'s  and de Finetti's theorem 
 that $\mu_N$ has weak infinite-volume limits (Gibbs measures)  which are convex 
 combinations of product measures \cite{SPOHN}. A variational characterization of these infinite-volume 
 measures and related results will be the content of another paper which will appear elsewhere. In 
 our current set-up we will always assume these infinite-volume measures exist for the class of 
 interactions we consider.

\subsection{Transforms of Mean-Field Models}\label{transforms}
We now introduce on $S\times S'$ a Borel probability measure $K$ such that 
\begin{equation}\begin{split}
K(d\xi_i)&=k(\s_i,\eta_i)\a(d\s_i)\a'(d\eta_i),\quad \text{with}\quad \sup_{(\s_i,\eta_i)\in 
S\times S'}\big|\log k(\s_i,\eta_i)\big|<\infty.
\end{split}
\end{equation}
We assume further that  
$ \a=\int k(\cdot,\eta_i)\a'(d\eta_i)$ and
$\a'=\int k(\s_i,\cdot)\a(d\s_i)$, where we are using the subscript  $i\in\N$ to convey the idea that 
the $K$ is the joint a priori measure for site $i$.
Given a Gibbs measure $\mu$ for the initial model \eqref{mean-field}, it is our aim in this work
 to investigate the Gibbs properties of the transformed measures
\begin{equation}\label{transmeasure} 
\mu'(d\eta)=\int_\O \mu(d\s)\prod_{i\in\N}k(\s_i,\eta_i)\a'(d\eta_i)
\end{equation}
as has been done  for the corresponding short-range models in \cite{KULOP}.
The study, as in \cite{KULOP}, will be based on investigating the properties of the finite-volume
conditional distributions  of the transformed system. This consists in studying  the 
infinite-volume
$N-$limits of the following finite-volume quantities:
 \begin{equation}\begin{split}\label{finvol}
 \mu'_{n,N}\big(d\eta_{V_n}\big|\eta_{V_N\ba V_n}\big)&=\int_{S^N}\mu_N(d\s_{V_N})\prod_{j=n+1}
 ^N k(\s_j,\eta_j)
 \prod_{i=1}^n k(\s_i,\eta_i)\a'(d\eta_i),
 \end{split}
 \end{equation}
 for fixed $n\in\N$ with $1\leq n<N$. But unlike in the lattice
spin systems the boundary conditions here will be fixed up to permutations, i.e., each 
boundary condition $\eta_{V_N\ba V_n}$ 
will be a representative of a class of configurations which gives rise to the same 
empirical measure.
In view of this, we shall take probability measures in $\PP(S')$ as our boundary 
conditions. As we shall show 
below, the infinite-volume $N$-limit of $\mu'_{n,N}$ will always factorize and this 
factorization necessitates the study of the $n=1$ case.  This leads to
the following definition of Gibbsianness for mean-field models which was originally 
 introduced by one of the
authors of this paper for the corresponding Curie-Weiss model \cite{KUL2,KUL1}. 
The case studied here
is a generalization of this notion from empirical average to empirical measures.

\begin{defn}\label{gibbs}
We call $\nu'\in\PP(S') $ a good configuation if and only if 
\begin{enumerate}
\item the limit
 \begin{eqnarray}
\g'_{1}(d\eta_1|\l'):=\lim_{N\uparrow\infty}\mu'_{1,N}\left(d\eta_1\big|\eta_
{V_N\ba\{1\}}\right)\quad\text{where}\quad \l'=\lim_{N\uparrow\infty}\dfrac{1}{N}\sum_
{i=2}^{N}\delta_{\eta_i}
\end{eqnarray} 
exists for all $\l'$ in a weak neighborhood of $\nu'$ and
\item for any Borel subset $A\subset S',$ the 
function $\l'\mapsto\g'_{1}(A|\l')$ is weakly continuous at $\l'=\nu'$.
\end{enumerate}
We say $\mu'$ is \textbf{Gibbs} iff every configuration is good.
\end{defn}
In what follows (unless otherwise stated) continuity of maps on $\PP(S')$ will always be w.r.t. the weak 
topology.

In our  investigation of the  continuity properties of the single-site kernels $\g'_1$ for the 
transformed system we employ the machinery of large deviations theory. In view of this, 
we will recall some basic facts about large deviations theory that we will need in our analysis in 
the next subsection.

\subsection{Some facts about Large Deviations Theory} 
In this subsection we recall some facts about large deviations theory and  for detailed
discussion on this theory and its application to statistical mechanics
we refer the reader to \cite{COM89,ELLIS85}.
Let $X$ be a Polish space equipped with its Borel $\s-$algebra. 
\begin{defn}\label{LDP}
A sequence of probability measures $(Q_N)_{N\in\N}$ in $\PP(X)$ is said to satisfy
a {\em large deviation principle}  (LDP) on $X$ with rate $a_N$ (sequence of positive 
numbers tending to infinity) and rate function $I:X\rightarrow [0,+\infty]$ if 
\begin{enumerate}
	\item $I$ is lower semi-continuous on $X$, and the level sets $\{x\in X:I(x)\leq a\}$   
	are compact for all $a\in [0,+\infty)$;
	\item for any Borel subset $B$ of $X$,
	\begin{equation}
	-I(\mathring{B})\leq \liminf_{N\rightarrow \infty}a_N^{-1}\log Q_N(B)\leq \limsup_
	{N\rightarrow \infty}a_N^{-1}\log Q_N(B)\leq -I(\bar B),
\end{equation}
where for any subset $C$ of $X$, $I(C)=\inf_{x\in C}I(x)$, and $\mathring{C}$  and $\bar C$ are 
respectively the interior and the closure of $C$. 
\end{enumerate}
\end{defn}
As an example take $(Y_n)_{n\in\N}$, an i.i.d. sequence of random variables on $X$ with  $\r$
 as the law of $Y_1$. Let $Q_N$ be  the distribution of the empirical measures
$L_N=\frac{1}{N}\sum_{i=1}^N\d_{Y_i}$. Then $Q_N$ satisfies LDP with rate $N$ and rate function
\begin{equation}
S(\nu|\r)=\left\{
\begin{array}{rl}
\int \frac{d\nu}{d\r}\log\frac{d\nu}{d\r}d\r & \text{if }\quad \nu \ll \r \quad \text{and}
\quad\frac{d\nu}{d\r}\log\frac{d\nu}{d\r}\in\LL^1(\r) \\\\
+\infty & \text{otherwise}.
\end{array} \right.
\end{equation}
The above example is Sanov's theorem in large deviations theory as can be found e.g. in
Theorem II.4.3 of \cite{ELLIS85}.

Another important fact about LDP that we shall employ in our study is the {\em contraction 
principle} (see e.g. Theorem II.5.1 of \cite{ELLIS85}), which comes to play when one is 
concerned with partial summary of the information
weighted by $Q_N$. More 
precisely, suppose $\psi$ is a continuous function from the Polish space $X$ to another
Polish space $Y$ and $Q_N$ is a sequence of probability measures on $X$ satisfying the LDP with 
rate $a_N$ rate function $I$. Then the sequence $\hat Q_N=Q_N\circ\psi^{-1}$ of probability 
measures on $Y$ also satisfies LDP with rate $a_N$ and rate function $\hat I$ given by
\begin{equation}\label{contraction}
	\hat I(y)=\inf\big\{I(x):\psi(x)=y\big\}.
\end{equation}

Our last fact from LDP  concerns  the integrals of  exponentials of
functionals of random variables whose distributions satisfy LDP. 
This is  found e.g. in \cite{ELLIS85} as  Theorem II.7.2a. The result in 
\cite{ELLIS85} is more general than what is stated here.

\begin{fact}\label{mainprop}
Let $X$ be a Polish space and $Q_N$ a sequence of probability measures on $X$ obeying LDP 
with rate $a_N$ and rate function $I$. Suppose  that $F:X\rightarrow\R$, which is  continuous and
bounded below, and for each $N\in\N$ the integral $\int_X \exp\big(-a_N F(x)\big)Q_N(dx)$ is finite. Let
  $Q_{N,F}$ be the sequence of probability measures  given by
\begin{equation}
	Q_{N,F}(A)=\dfrac{\int_A \exp\big(-a_N F(x)\big)Q_N(dx)}{\int_X \exp\big(-a_N F(x)\big)Q_N(dx)},
\end{equation}
for any Borel subset $A$ of $X$. Then $Q_{N,F}$ satisfies LDP with rate $a_N$ and rate function
\begin{equation}\label{rateF}
	I_F(x)= I(x)+F(x)-\inf_{y\in X}\big(I(y)+F(y)\big).
\end{equation}
 \end{fact}

\section{Two-layer Models and Gibbsianness of Transformed Systems}\label{s3}
In this section we consider  mean-field systems with $S\times S'$ as their single-site spin space
and $K$ as the corresponding a priori measure. As is discussed in \cite{KULOP} the a priori measure
$K$ couples two systems; namely, the first-layer system described by some given mean-field interaction 
$\P$ an the a priori measure $\a$, and the second-layer system which are i.i.d. with distribution $\a'$. 
The finite-volume Gibbs measures $\td\mu_N$ for our two-layer (joint) system are given by
\begin{equation}\begin{split}\label{two-layer}
\td\mu_N(d\xi_{V_N})&=\dfrac{\exp\Bigl(-N\P\bigl(L_N(\s_{V_N})\bigr)\Bigr)\prod_{i=1}^NK(d\s_i,d\eta_i)}
{\int_{S^N}\exp\Bigl(-N\P\bigl(L_N(\hat\s_{V_N})\bigr)\Bigr)\prod_{i=1}^N\a(d\hat\s_i)}\cr
&=\dfrac{\exp\Bigl(-N\Bigl\{\P\bigl(\pi_1 L_N(\xi_{V_N})\bigr)-L_N(\xi_{V_N})\bigl[
\log k(\cdot,\cdot)\bigr]\Bigr\}\Bigr)\prod_{i=1}^N\a(d\s_i)\a'(d\eta_i)}
{\int\exp\Bigl(-N\Bigl\{\P\bigl(\pi_1 L_N(\hat{\xi}_{V_N})\bigr)-L_N(\hat{\xi}_{V_N})\bigl[
\log k(\cdot,\cdot)\bigr]\Bigr\}\Bigr)\prod_{i=1}^N\a(d\hat\s_i)\a'(\hat\eta_i)},
\end{split}
\end{equation}
 where $L_N(\xi_{V_N})=\frac{1}{N}\sum_{i=1}^N\d_{(\s_i,\eta_i)}$ is the 
 joint empirical measures
 and $\pi_1 L_N(\xi_{V_N})=L_N(\s_{V_N})$ is the projection onto the first 
 variable. We have also 
 denoted by $\nu[f]$ the integral of the measurable map $f$ w.r.t. the measure 
 $\nu$. Under  the map $\xi_{V_N}\mapsto L_N(\xi_{V_N})$, the joint 
 measures $\td\mu_N$ have unique
 push-forwards $\td Q_{N,\P,k}$ in $\PP\bigl(\PP(S\times S')\bigr)$. So we now study 
the large deviation properties  of the empirical measures $L_N(\xi_{V_N})$  under 
$\td Q_{N,\P,k}$ and the proposition below is a summary of this LDP property.
 
 \begin{prop}\label{jointldp}
 The sequence of probability measures $\td Q_{N,\P,k}$ satisfies an LDP with rate 
 $N$ and rate function $\td J$ given by
\begin{equation}\begin{split}\label{jointrf}
	\td J(\td\nu)&= \P(\pi_1\td\nu)-\td\nu\bigl[\log k(\cdot,\cdot)\bigr]+ S(\td\nu|
	\a\otimes\a')-{\rm const}\quad \text{where}\cr
{\rm const}&=\inf_{\td\l\in\PP(S\times S')}\Bigl\{\P(\pi_1\td\l)-\td\l\bigl[\log 
k(\cdot,\cdot)\bigr] +S(\td\l|\a\otimes\a') \Bigr\}.
\end{split}
\end{equation}
 
 \end{prop}
\textbf{Remark:} 
\textbf{1. } Proposition \ref{jointldp} is a direct consequence of Fact \ref{mainprop}, 
since $\P$ is bounded and  the weak topology on $\PP(S\times S')$ turns $\PP(S\times S')$ into 
a compact separable 
metric space.\\
\textbf{2. } The continuous map  $\xi_{V_N}\mapsto L_N(\xi_{V_N})$ from 
$(S\times S')^N$ to $\PP(S\times S')$ also 
induces a continuous map (say $\vartheta_N$) from $\PP\big((S\times S')^N\big)$ into $\PP
\bigl(\PP(S\times S')\bigl)$. So the push-forwards
$\td Q_{N,\P,k}$ are simply obtained by replacing  the product measures $\prod_
{i=1}^N\a\otimes\a'$ in the second 
equality in \eqref{two-layer}  by their images under 
$\vartheta_N$.\\
\textbf{3. } Observe from Proposition \ref{jointldp} that the empirical measures  of the 
initial spin variables under the initial finite-volume Gibbs measures $\mu_N$ \eqref{mean-field} satisfy an
LDP with rate $N$ and rate function $I_{\a}^\P$ given by
\begin{equation}\label{initialrf}
I^\P_\a(\nu)= S(\nu|\a)+\P(\nu)-\inf_{\mu\in\PP(S)}\Big[S(\mu|\a)+\P(\mu)\Big].
\end{equation}

Our next task is to study the LDP for the transformed measures $\mu'_N$ given by
\begin{equation}
	\mu'_N=\int_{S^N}\td\mu_N(d\s_{V_N}).
\end{equation}
For any $\td\nu\in\PP(S\times S')$ we denote by $\pi_2\td\nu$  the marginal of $\td\nu$
on $S'$. Denote  by $Q'_{N,\P,k}$  the push-forward of $\mu'_N$ given by  
$Q'_{N,\P,k}=\td Q_{N,\P,k}\circ \pi_2^{-1}$. Therefore by the contraction principle, 
the empirical measures of the transformed spins under $Q'_{N,\P,k}$  satisfy an LDP which 
we formulate in the proposition below.
\begin{prop}\label{second-layerldp}
 The sequence of probability measures $Q'_{N,\P,k}$ in $\PP\bigl(\PP(S')\bigr)$ satisfies an LDP
 with rate $N$ and rate function $J'$ given by 
 \begin{equation}\begin{split}\label{slayerrf}
 J'(\nu')&=\inf_{\substack{\td\nu\in\PP(S\times S')\\ \pi_2\td\nu=\nu'}} \td J(\td\nu)=S(\nu'|\a')+
 \inf_{\td\nu\in M_{\nu'}}J_{\nu'}(\td\nu)-{\rm const},\quad \text{where}\cr
 	J_{\nu'}(\td\nu)&=S\big(\td\nu|\a\otimes\nu'\big)+\P(\pi_1\td\nu)-
	\td\nu\bigl[\log k(\cdot,\cdot)\bigr],
\end{split}
 \end{equation}
and  $M_{\nu'}$ is the subset of $\PP(S\times S')$ consisting of probability measures  with 
 fixed second marginal $\nu'$.
 \end{prop}

\begin{remark}
For  each $\nu'\in\PP(S')$, $J_{\nu'}$ up to an additive constant (depending on 
$\nu'$) is the large deviations
rate function  for the joint system when the second-layer system is constrained  to configurations 
with empirical measure $\nu'$, i.e. $J_{\cdot}$ (up to an additive constant)  is the rate function
for the CFLM. 
This CFLM rate function will play a key role in determining whether  
the transformed system is Gibbs  or not. We  shall show below that the continuity properties 
of the transformed single-site kernels $\g'_1(\cdot|\nu')$ will be determined by the unicity 
of the global minimizers of the  function $J_{\nu'}$ uniformly in 
$\nu'$. Observe also that $J_{\nu'}$ is a lower semi-continuous function and that it  also attains 
its infimum on  $M_{\nu'}$, since $M_{\nu'}$ is compact subset of $\PP(S\times S')$. Additionally, 
$M_{\nu'}$ is  convex. 
\end{remark}
\textbf{Proof of Proposition \ref{second-layerldp}: } 
The first equality in the expression for $J'$  \eqref{slayerrf} of  Proposition 
\ref{second-layerldp} follows from  the contraction principle, since the 
map $\pi_2:\PP(S\times S')\rightarrow\PP(S')$ is weakly  continuous. Further, for each $\nu'\in\PP(S')$ 
the set $M_{\nu'}$ is compact, because  $M_{\nu'}=\pi_2^{-1}\bigl(\{\nu'\}\bigr)$ and by 
continuity of $\pi_2$ it is a closed subset of $\PP(S\times S')$. Now it also follows from standard 
results in analysis that closed subsets of compact set are compact and hence  the compactness of $M_{\nu'}$, 
since the weak topology turns $\PP(S\times S')$ into a compact Polish space by the compactness 
of  $S\times S'$. This shows that the infimum of $\td J$ over measures with fixed second marginal $\nu'$ is
attained on $M_{\nu'}$.

Now for each $\;\nu'\in\PP(S')$, the measures 
$\;\td\nu \in M_{\nu'}\;$ are of the form 
$\td\nu(d\xi_i)=\nu'(d\eta_i)\td\nu(d\s_i|\eta_i)$, where for each $\eta_i\in S'$, $\;\td\nu
(d\s_i|\eta_i)\in \PP(S)$. With this representation, the relative entropy for 
the elements in $M_{\nu'}$ w.r.t. $\a\otimes\a'$ takes the form
\begin{equation}\begin{split}
	S(\td\nu|\a\otimes\a')&=S(\nu'|\a') +\int_{S'}\nu'(d\s')S\bigl(\td\nu(\cdot|\s')\big|\a'\bigr)
	=S(\nu'|\a')+ S(\td\nu|\a\otimes \nu').
\end{split}
\end{equation}
This proves  the second equality of the expression for $J'$ in \eqref{slayerrf}.
\begin{flushright}
$\Cox$
\end{flushright}

 The measures  in $M_{\nu'}$ that are of interest in determining $J'(\nu')$  are at
 most those  for which $\td\nu(\cdot|\eta_i)\ll \a$, for all $\eta_i\in S'$ i.e.,
\begin{equation}\label{18nu}
	\td\nu(d\xi_i)=\nu'(d\eta_i)\a(d\s_i)f_{\nu'}(\s_i|\eta_i),
\end{equation}
for some measurable function $f_{\nu'}:S\times S'\rightarrow [0,\infty]$ with the property that
for each $\eta_i\in S', \quad \int \a(d\s_i) f_{\nu'}(\s_i|\eta_i)=1 $. We call $f_{\nu'}$ the 
{\em conditional $\a$- density} of $\td\nu$.
This reduces the whole problem of minimizing $J_{\nu'}$ over $M_{\nu'}$ to the problem of finding the 
conditional $\a-$densities $f_{\nu'}$ for which $\nu'(d\eta_i)\a(d\s_i)f_{\nu'}(\s_i|\eta_i)$ is a 
minimizer of $J_{\nu'}$.

As our next result, we present an explicit expression  for the conditional $\a-$densities 
$f_{\nu'}(\s_i|\eta_i)$ at 
which $J_{\nu'}$ attains  both global and local minima. 
\begin{thm}\label{minimizers}
For any $\nu'\in\PP(S')$, the function $J_{\nu'}$ attains its infimum on $M_{\nu'}$. Furthermore, 
any  minimizer 
(global or local) $\td\nu\in M_{\nu'}$ of  $J_{\nu'}$   has the conditional $\a$-density
 $f_{\nu'}$  which is  $\;\a\otimes\nu'$-a.s. strictly positive and  satisfies the 
"constrained mean-field equation" 
\begin{eqnarray} \label{cdf}
f_{\nu'}(\s_i|\eta_i)=\dfrac{e^{-\P^{(1)}(\pi_1\td\nu,\d_{\s_i})}k(\s_i,\eta_i)}{\int e^{-\P^{(1)}(\pi_1\td\nu,\d_{\hat\s_i})}k(\hat\s_i,\eta_i)\a(d\hat{\s_i})}.
\end{eqnarray}
 \end{thm}

\begin{remark}\textbf{1. } As we pointed out in the above, the measures $\td\nu\in M_{\nu'}$ that are involved in
determining $J'(\nu')$ are those that takes the form \eqref{18nu}. The minimizers of $J_{\nu'}$ are 
among these probability measures and indeed they are  those probability measures with $f_{\nu'}$ 
given by \eqref{cdf}.\\
\textbf{2. } Note that the minimizers of the rate function $I_\a^\P$ \eqref{initialrf} for the initial system are also measures
$\nu\in\PP(S)$ with $\nu(d\s_i)=f(\s_i)\a(d\s_i)$, where the $f$'s are $\a$-a.s. strictly positive and statisfy the
mean-field equation 
\begin{eqnarray} \label{density}
f(\s_i)=\dfrac{e^{-\P^{(1)}(\nu,\d_{\s_i})}}{\int e^{-\P^{(1)}(\nu,\d_{\hat\s_i})}\a(d\hat{\s_i})}.
\end{eqnarray}
\end{remark}
We defer the proof of Theorem \ref{minimizers} to the Appendix at the end of the paper.

The above  consideration leads to the following definition of some probability kernel from $M_{\nu'}$ to  
$S\times S'$.
\begin{defn}\label{cflpk}
We refer to the map $\g_{\nu'}:M_{\nu'}\rightarrow M_{\nu'}$  given by
\begin{equation}\label{rckernel}
\g_{\nu'}(d\xi_i|\td\nu)=\nu'(d\eta_i)\dfrac{e^{-\P^{(1)}(\pi_1\td\nu,\d_{\s_i})}k(\s_i,\eta_i)\a
(d\s_i)}{\int e^{-\P^{(1)}(\pi_1\td\nu,\d_{\hat\s_i})}k(\hat\s_i,\eta_i)\a(d\hat{\s_i})}
\end{equation}
as the {\em constrained first layer probability kernel} (CFLPK).
\end{defn}
Note that for any Borel subset $A$ of $S\times S'$, the map $\g_{\nu'}(A|\cdot):M_{\nu'}
\rightarrow \R$ is 
continuous by the continuity property of $\P^{(1)}$. Observe further that  not 
all the  measures $\td\nu\in M_{\nu'}$   have $\td\nu(\cdot|\cdot)\ll\a$, but the fixed points of 
$\g_{\nu'}$ do and include all the minimizers of $J_{\nu'}$.  

It is not hard to deduce from the second remark below Theorem \ref{minimizers} that we can analogously 
define a version of the CFLPK for the initial system, namely;
\begin{equation}\label{kernel}
\g_1(d\s_i|\nu)=\frac{e^{-\P^{(1)}(\nu,\d_{\s_i})}\a
(d\s_i)}{\int e^{-\P^{(1)}(\nu,\d_{\hat\s_i})}\a(d\hat{\s_i})}.
\end{equation}
Here $\g_1$ is a map from $\PP(S)$ onto itself.

We introduce the following notion   of consistency for mean-field models.
\begin{defn}\label{consistency}
A probability measure $\td\nu\in M_{\nu'}$ is said to be consistent w.r.t. $\g_{\nu'}$ 
whenever
\begin{equation}
\g_{\nu'}(\xi_i|\td\nu)=\td\nu(d\xi_i).
\end{equation}
\end{defn}
The consistent probability measures for the CFLPK are those measures $\td\nu\in M_{\nu'}$  for which 
$\td\nu(\cdot|\cdot)\ll\a$ and have conditional $\a-$density functions $f_{\nu'}$ which are given 
by \eqref{cdf}. 
Observe that the consistent probability measures are the fixed points of $\g_{\nu'}$. 

This notion of consistency also carries over to the kernel $\g_1$ for the initial system. Observe further  
that the kernel $\g_1(\cdot|\nu)$ is continuous as a function of the conditioning $\nu\in\PP(S)$. 
This continuity property  follows as a result of the continuity property of $\P^{(1)}(\cdot,\cdot)$  in its first 
argument. This then implies the Gibbsianness of the initial system by Definition \ref{gibbs}.

\begin{prop}\label{tikhonov}
For any probability measure $\nu'\in \PP(S')$,  the CFLPK $\g_{\nu'}$ has a fixed point. 
\end{prop}
\textbf{Proof:} The existence of a fixed point for the CFLPK follows from the Tychonov's fixed point 
theorem,  which states that for any non-empty compact convex subset $X$ of a locally convex 
topological vector space $V$, and continuous function $f:X \rightarrow X$, there is a fixed point for f.

Now as we observed from the proof of Proposition \ref{second-layerldp}, for any $\nu'\in \PP(S')$ the set
$M_{\nu'}$ is compact w.r.t. the weak topology. It is also not hard to see that it is  convex. So we are
now only left to show that the space $\MM(S\times S')$ of finite signed measures on $S\times S'$ is locally 
convex topological vector space under the weak topology.

{\em That} $\MM(S\times S')$ is a locally convex topological vector space w.r.t. the weak topology follows from the
following sequence of arguments:
The total variational norm  turns the space $\MM(S\times S')$ into a Banach space. This then implies that $\MM(S\times S')$
is normable and consequently it is  locally convex (i.e. the  origin has a local base of convex sets)  w.r.t. the 
total variational topology by Theorem 1.39 of \cite{RUD}.  It then  follows from the corollary to Theorem 3.4 of  
\cite{RUD} that the dual space  $\MM^*(S\times S')$ of $\MM(S\times S')$ separates points in $\MM(S\times S')$.
The weak topology on $\MM(S\times S')$ generated by the dual space $\MM^*(S\times S')$ then   turns 
$\MM(S\times S')$ into a locally convex topological space by Theorem 3.10 of \cite{RUD}.

\begin{flushright}
$\Cox$
\end{flushright}
The existence of consistent probability measures for the initial kernel $\g_1$ is trivial, since $\g_1$ is continuous
and $\PP(S)$ is compact convex set. Note also that for each $\nu'\in\PP(S')$, the constrained first layer model is Gibbsian (in the sense of Definition \ref{gibbs}) by the continuity property of $\g_{\nu'}$. Indeed, the CFLPK $\g_{\nu'}$ can be extended to the whole of $\PP(S\times S')$ but those probability measures we will be interested in 
are those in $M_{\nu'}$.

For each $\nu'\in \PP(S')$ denote by ${\cal C}_{\nu'}$ the set of all consistent probability 
measures   of $\g_{\nu'}$. We state as our next result the following  lemma  concerning a single-site
 variational principle for the CFLM. 
\begin{lem}\label{rcentropy}
For any probability measure $\td\nu\in M_{\nu'}$ the relative entropy $S(\td\nu|\a\otimes \nu')$ 
satisfies
\begin{equation}\label{vprin}
S(\td\nu|\a\otimes \nu')\geq -\P^{(1)}(\pi_1\td\nu,\pi_1\td\nu)+\td\nu\left[\log k(\cdot,\cdot
)\right]-\int\nu'(d\eta_i)\log\int\a(d\hat{\s_i}) e^{-\P^{(1)}(\pi_1\td\nu,\d_{\hat\s_i})}k(\hat\s_i,\eta_i).
\end{equation}
In particular equality is attained whenever $\td\nu\in {\cal C}_{\nu'}$.
\end{lem}
\textbf{Proof:}
For any $\td\nu\in M_{\nu'}$ the expression on the right-hand-side of \eqref{vprin} becomes
\begin{equation}\begin{split}
 &-\P^{(1)}(\pi_1\td\nu,\pi_1\td\nu)+\td\nu\left[\log k(\cdot,\cdot
)\right]-\int\nu'(d\eta_i)\log\int \a(d\s_i) e^{-\P^{(1)}(\pi_1\td\nu,\d_{\s_i})}k(\s_i,\eta_i)\cr
&= \int_{S\times S'} \td\nu(d\xi_i)\log\left\{\frac{e^{-\P^{(1)}(\pi_1\td\nu,\d_{\s_i})}k(\s_i,\eta_i)}
{\int \a(d\s_i) e^{-\P^{(1)}(\pi_1\td\nu,\d_{\s_i})}k(\s_i,\eta_i)}\right\}=\int_{S\times S'} 
\td\nu(d\xi_i)\log\left\{\frac{d\g_{\nu'}}{d(\a\otimes\nu')}(\xi_i|\td\nu)\right\}.
\end{split}
\end{equation}
The proof  now becomes showing that 
\begin{equation}\begin{split}
S(\td\nu|\a\otimes \nu')- \int_{S\times S'} 
\td\nu(d\xi_i)\log\left\{\frac{d\g_{\nu'}}{d(\a\otimes\nu')}(\xi_i|\td\nu)\right\}\geq0.
\end{split}
\end{equation}
The case for measures $\td\nu\in M_{\nu'}$  with  $S(\td\nu|\a\otimes \nu')=\infty$ is trivial
since we get strict inequality by the  boundedness properties of $\P$ and $k$.

Now for the case  of $\td\nu\in M_{\nu'}$  with  $S(\td\nu|\a\otimes \nu')<\infty$, we obtain
\begin{equation}\begin{split}
&\int_{S\times S'} 
\td\nu(d\xi_i)\log\left\{\frac{d\td\nu}{d(\a\otimes\nu')}(\xi_i)\right\}- \int_{S\times S'} 
\td\nu(d\xi_i)\log\left\{\frac{d\g_{\nu'}}{d(\a\otimes\nu')}(\xi_i|\td\nu)\right\}\cr
& =\int_{S\times S'} \td\nu(d\xi_i)\log\left\{\frac{d\td\nu}{d\g_{\nu'}(\cdot|\td\nu)}(\xi_i)\right\}
=S(\td\nu\big|\g_{\nu'}(\cdot|\td\nu))\geq 0,
\end{split}
\end{equation}
 since $S(\td\nu\big|\g_{\nu'}(\cdot|\td\nu))\geq 0$  with equality holding only when $\td\nu=
 \g_{\nu'}(\cdot|\td\nu)$. This concludes the proof.
\begin{flushright}
$\Cox$
\end{flushright}

We state as our next result a theorem concerning some function $\Psi_{\nu'}$ on 
$M_{\nu'}$ which is dominated by   $J_{\nu'}$ and coincides with $J_{\nu'}$ on 
${\cal C}_{\nu'}$.

\begin{thm}\label{dom}
For any given $\nu'\in \PP(S')$, the  function $J_{\nu'}$ satisfies
\begin{equation}\begin{split}\label{jmin}
J_{\nu'}(\td\nu)&\geq \Psi_{\nu'}(\td\nu)\quad \text{where}\cr
\Psi_{\nu'}(\td\nu)&=\P(\pi_1\td\nu)-\P^{(1)}(\pi_1\td\nu,\pi_1\td\nu)
-\int\nu'(d\eta_i)\log\int \a(d\hat{\s_i}) e^{-\P^{(1)}(\pi_1\td\nu,\d_{\hat\s_i})}k(\hat\s_i,\eta_i).
\end{split}
\end{equation}
In particular $J_{\nu'}$ coincides with $\Psi_{\nu'}$ on ${\cal C}_{\nu'}$ and  if $\P$ is 
homogeneous of degree $p$ then $\Psi_{\nu'}$  becomes
\begin{equation}\begin{split}\label{jmin1}
\Psi_{\nu'}(\td\nu)&=(1-p)\P(\pi_1\td\nu)-\int\nu'(d\eta_i)\log\int\a(d\hat{\s_i}) e^{-\P^{(1)}(\pi_1\td\nu,\d_
{\hat\s_i})}k(\hat\s_i,\eta_i).
\end{split}
\end{equation}
\end{thm}
\textbf{Proof:}
The expression for $\Psi_{\nu'}$ and the inequality \eqref{jmin} follow by substituting the 
 lower bound on $S(\cdot|\a\otimes\nu')$ in Lemma \ref{rcentropy} into the expression for $J_{\nu'}$ in 
 Proposition \ref{second-layerldp}. 
Furthermore, if $\P$ is homogeneous of degree $p$,  we then have  
\begin{equation}
\P^{(1)}(\pi_1\td\nu,\pi_1\td\nu)=\frac{d}{dt}\P(\pi_1\td\nu +t\pi_1\td\nu)\Big|_{t=0}=
\frac{d}{dt}(1+t)^p\P(\pi_1\td\nu )\Big|_{t=0}=p\P(\pi_1\td\nu )
\end{equation}
and putting this into the expression for $\Psi_{\nu'}$ in \eqref{jmin} yields the desired expression 
in \eqref{jmin1}.
\begin{flushright}
$\Cox$
\end{flushright}

\begin{cor}
The transformed 
LDP rate function $J'$ now becomes
\begin{equation}\begin{split}\label{transrf1}
J'(\nu')&=S(\nu'|\a')+\P_{k}(\nu')-{\rm const},\quad \text{where}\quad
\P_{k}(\nu')=\inf_{\td\nu\in{\cal C}_{\nu'}}\Psi_{\nu'}(\td\nu).
\end{split}
\end{equation} 
\end{cor}
\begin{remark}
\textbf{1. } The above expression for the transformed rate function $J'$ is a consequence of
the fact that  $J_{\nu'}$ coincides with $\Psi_{\nu'}$ on ${\cal C}_{\nu'}$. 
Thus $\P_k$ is the interaction for the transformed system arising from the initial system 
described by $\P$ and subjected to the site-wise transformations governed by $k$.
\end{remark}

\subsubsection{Examples}
Take $\P$ to be  an Ising  mean-field interaction (i.e. $S=\{+1,-1\}$) given by 
\begin{equation}
\P(m)=-\frac{\beta}{p}m^p
\end{equation}
where $m\in[-1,1]$  and  $p\geq1$. Here the reason for  using  $m$  instead of probability measures 
on $S$ is that the probability
measures on $S$ are uniquely determined by $m$, i.e. each $m\in[-1,1]$ can uniquely be associated 
with a probability measure (say $\nu$) on $S$ given by $\nu(\s_i)=\frac{1+m}{2}
\d_{+1}(\s_i)+\frac{1-m}{2}\d_{-1}(\s_i)$. The 
expectation w.r.t. $\nu$ then gives rise to $m$.

We take  $k(\s_i,\eta_i)
=p_t(\s_i,\eta_i)$ to be  the transition probabilities (i.e. $p_t(\s_i,\eta_i)$ is the  
probability of starting with $\s_i$ at site $i$ and observing $\eta_i$  after $t$ time units) for rate one site-wise 
independent 
spin-flip dynamics on $S$  \cite{KUL2}.  Here both $S$ and $S'$ are the same, and the a priori measures $\a=\a'
=\frac{1}{2}(\d_{+1}+\d_{-1})$.  
More precisely, $p_t(\s_i,\eta_i)$
is given by 
\begin{equation}
p_t(\s_i,\eta_i)=\frac{e^{\s_i\eta_i h_t}}{2 \cosh(h_t)},\quad \text{where}\quad h_t=\frac{1}{2}\log\frac{1-
e^{-2t}}{1+e^{-2t}}.
\end{equation}
As pointed out above, we will denote by $\tau\in[-1,1]$ the expected values of the  probability measures on the
 transformed single-site space. We will also write $m'$ for the expected values of the first marginals
of probability measures on $S\times S$ with fixed second marginal. The fixed second  marginal will be assumed 
to have mean $\tau$.

Then $\Psi_{\cdot}$ for this set-up becomes
\begin{equation}\begin{split}
\Psi_{\tau}(m')=&\frac{(p-1)\beta}{p} {m'}^p- 
\frac{1+\tau}{2}\log\left\{\cosh\left(\beta {m'}^{(p-1)}+h_t\right)\right\}\cr-&\frac{1-\tau}{2}\log\left\{\cosh\left(\beta {m'}^{(p-1)}-h_t\right)\right\}+\log\left(2\cosh(h_t)\right).
\end{split}
\end{equation}
Consequently, this form of $\Psi_\tau$ gives rise to the mean-field equation
\begin{equation}
m'= \frac{1+\tau}{2}\tanh\left(\beta {m'}^{(p-1)}+h_t\right)+\frac{1-\tau}{2}\tanh\left(\beta {m'}^{(p-1)}-h_t\right).
\end{equation}

\begin{remark} 
In the case $p=2$, $\Psi_\cdot$ is the Hubbard-Stratonovitch  potential function \cite{KUL2}. 
The unicity of global minimizers
of this potential function played a crucial role in  determining the Gibbs and non-Gibbs 
properties of the corresponding  transformed system studied in \cite{KUL2}.

In \cite{KUL2} the derivation of $\Psi_\cdot$  is based on the quadratic nature of the 
interaction $\P$. The technique employed there cannot be used to derive $\Psi_\cdot$ for non-quadratic 
interactions and this is where our approach comes to the rescue, i.e. our approach of deriving $\Psi_
{\cdot}$ via the machinery of large deviations is adaptable to more
general mean-field interactions (both discrete and continuous spins). 
 \end{remark}
 We now discuss  in detail Gibbsianness and non-Gibbsianness for mean-field models as introduced in
 Definition  \ref{gibbs}.
 \subsection{Gibbsianness for Transforms of mean-field models}
 In this subsection we study the Gibbs properties of the transformed measures 
 $\mu'$ \eqref{transmeasure} introduced in Section \ref{transforms}. This investigation of
 the Gibbs properties of transformed measures shall be based on the continuity properties of the
 conditional distributions $\g'$ as a function of the conditioning. Before we formulate this, 
 let us fix some notations that we shall
 use in our formulation. For each $N\geq 2$  and $1\leq n<N $ we denote by 
 $\td\mu_{N- n}[\bar\eta_{V_N\ba V_n}]$  the joint system in $V_N\ba V_n$ when  the second-layer spins 
 are constrained to a given configuration $\bar\eta\in \O'$, i.e. $\bar\eta_{V_N\ba V_n}$ is the projection 
 of $\bar\eta$ onto ${S'}^{N-n}$. As we pointed out in the above,  a representative $\bar\eta_{V_N\ba V_n}$
 of a class of configurations in ${S'}^{N-n}$ with the same empirical measure $L_{N-n}\big(\bar\eta_
 {V_N\ba V_n}\big)$ will give rise to the same measure $\td\mu_{N- n}[\bar\eta_{V_N\ba V_n}]$. 
 Therefore by fixing $\bar\eta_{V_N\ba V_n}$  implies we are restricting attention
 to only the configurations in a subset of ${S'}^{N-n}$ with  fixed $L_{N-n}\big(\eta_{V_N\ba V_n}\big)$. 
 Suppose $\bar\eta_{V_N\ba V_n}$ is one such representative, we call $\td\mu_{N- n}[\bar\eta_{V_N\ba V_n}]$  the
{\em restricted constrained first layer model (RCFLM)} for the corresponding mean-field model. 
It is is restricted because we are not taking into account the spins in $V_n$ and constrained because we have frozen the
configurations in the second-layer to $\bar\eta_{V_N\ba V_n}$. 
 More precisely, 
 
\begin{equation}\begin{split}\label{rcflm}
	\td\mu_{N-n}[\bar\eta_{V_N\ba V_n}](d\s_{V_N\ba V_n})&=\dfrac{\exp\Bigl( -N\P\bigl(\pi_1\bar\nu_{N,n}\bigr)
	\Bigr)\prod_{i=n+1}^N k(\s_i,\bar\eta_i)\a(d\s_i)}{\int\exp\Bigl( -N\P\bigl(\pi_1\hat{\bar\nu}_{N,n}\bigr)
	\Bigr)\prod_{i=n+1}^N k(\hat\s_i,\bar\eta_i)\a(d\hat{\s}_i)}\cr
	&=\dfrac{\exp\Bigl( -N\Bigl\{\P\bigl(\pi_1\bar\nu_{N,n}\bigr)
	-\bar\nu_{N,n}
	\big[\log k\big]
	\Bigr\}\Bigr)\prod_{i=n+1}^N \a(d\s_i)}{\int\exp\Bigl( -N\Bigl\{\P\bigl(\pi_1\hat{\bar\nu}_{N,n}\bigr)
	-\hat{\bar\nu}_{N,n}\bigl[\log k\bigr]\Bigr\}\Bigr)\prod_
	{i=n+1}^N\a(d\hat{\s}_i)},
	\end{split}
\end{equation}
 where $\bar\nu_{N,n}=\frac{N-n}{N}L_{N-n}(\bar\xi_{V_N\ba V_n})$,
 $L_{N-n}(\bar\xi_{V_N\ba V_n})=\frac{1}{N-n}\sum_{i=n+1}^N\d_{(\s_i,\bar\eta_i)}
 $ and $\hat{\bar\xi}_i=(\hat\s_i,\bar\eta_i)$.   \\\\
\textbf{Remark:} Suppose $\nu'\in\PP(S')$ is the  empirical measure for the configuration $\bar\eta$. 
Then for a fixed $n$, the sequence of measures $\bar\nu_{N,n}$ under the pushed-forwards of  
$\td\mu_{N-n}[\bar\eta_{V_N\ba V_n}]$ satisfies an LDP with rate
$N$ and rate function $\bar J_{\nu'}:M_{\nu'}\rightarrow \R\cup \{+\infty\}$ given by   
\begin{equation}
\bar J_{\nu'}(\td\nu)=J_{\nu'}(\td\nu)-{\rm const}_{\nu'},\quad \text{where}\quad {\rm const}_{\nu'}=
\inf_{\td\nu\in M_{\nu'}}J_{\nu'}(\td\nu).
\end{equation}
The validity of the expression for $\bar J_{\nu'}$ lies in the fact that for $\bar\eta\in\O'$ constrained 
to have empirical measure $\nu'$,  the sequence of empirical measures $L_{N-n}(\bar\xi_{V_N\ba V_n})$ converge weakly
 in the $N$-limit to  an element in $M_{\nu'}$. For any  of such measures $\td\nu$ for which $\td\nu(\cdot|\eta_i)\ll\a$ 
(for $\nu'$ almost all $\eta_i\in S'$),
 \begin{equation}
 \frac{d\td\nu}{d\a}(\s_i,d\eta_i)=\nu'(d\eta_i)f_{\nu'}(\s_i|\eta_i). 
 \end{equation}  
 Therefore, the relative entropy of such probability measures  $\td\nu$ w.r.t. $\a$ then  becomes 
 \begin{equation}
 S(\td\nu|\a)=\int \nu'(d\eta_i)\int f_{\nu'_i}(\s_i|\eta_i)\log f_{\nu'_i}(\s_i|\eta_i)\a(d\s_i).
 \end{equation} 
 
Our next result in this subsection concerns a representation of the finite-volume transformed conditional 
distributions $\mu'_{n,N}(\cdot|\bar\eta_{V_N\ba V_n})$ in terms of the RCFLM  $\td\mu_{N-n}[\bar\eta
_{V_N\ba V_n}]$.

\begin{lem}\label{mainlemma}
Let $N,\; n \; \text{and}\; \bar\eta_{V_N\ba V_n}$ be as above. Then the finite-volume conditional 
distribution
$\mu'_{n,N}(\cdot|\bar\eta_{V_N\ba V_n})$ for the transformed system has the form
\begin{equation}
\mu'_{n,N}(d\eta_{V_n}|\bar\eta_{V_N\ba V_n})=\dfrac{\displaystyle{\td\mu_{N-n}[\bar\eta_{V_N\ba V_n}]\left
[\prod_{i=1}^n \int_{S} e^{-\P^{(1)}\left(\pi_1\bar\nu_{N,n},\d_{\s_i}\right)+o\big(\frac{1}{N}\big)}
k(\s_i,\eta_i)
\a(d\s_i)\a'(d\eta_i)\right]}}{\displaystyle{\td\mu_{N-n}[\bar\eta_{V_N\ba V_n}]\left[\prod_{i=1}^n \int_{S} e^{-\P^{(1)}\left(\pi_1\bar\nu_{N,n},\d_{\s_i}\right)+o\big(\frac{1}{N}\big)}
\a(d\s_i)\right]}}.
\end{equation}

\end{lem}
\textbf{Proof:}
Note from the definition of the the transformed system that we can write $\mu'_{n,N}(d\eta_{V_n}|\bar\eta_{V_N\ba V_n})$
as
\begin{equation}\label{defin}
\mu'_{n,N}(d\eta_{V_n}|\bar\eta_{V_N\ba V_n})=\dfrac{\displaystyle{\prod_{j=n+1}^N\int_{S} \a(d\s_j)k(\bar\xi_j) 
\prod_{i=1}^n \int_{S}\a(d\s_i)  \a'(d\eta_i)e^{-N\P\left(\pi_1 L_N(\bar\xi_{V_N})\right)}
k(\xi_i)
}}{\displaystyle{\prod_{j=n+1}^N \int_S\a(d\s_j)k(\bar\xi_j) \prod_{i=1}^n \int_{S\times S'}\a(d\s_i)\a'
(d\hat\eta_i) e^{-N\P\left(\pi_1 L_N(\bar\xi_{V_N})\right)}
k(\hat\xi_i)}},
\end{equation}
where the joint configuration $\bar\xi_{V_N}$ is such that $\bar\xi_{V_N\ba V_n}=(\s_i,\bar\eta_i)_
{i\in V_N\ba V_n }$ and 
$\bar\xi_{ V_n}=(\s_i,\eta_i)_{i\in  V_n }$, and $\hat\xi_i=(\s_i,\hat\eta_i)$. Now by writing the 
joint empirical measure as  $$L_N(\bar\xi_{V_N})=
\frac{N-n}{N}L_{N-n}(\bar\xi_{V_N\ba V_n})+\frac{n}{N}L_n(\xi_{V_n}),$$ and adding and subtracting 
$N\P(\frac{N-n}{N}L_{N-n}(\bar\xi_{V_N\ba V_n}))$ from the exponent $N\P(L_{N}(\bar\xi_{V_N}))$ we obtain
\begin{equation}\begin{split}\label{decom}
N\P\left(L_{N}(\bar\xi_{V_N})\right)=&N\P\left(\frac{N-n}{N}L_{N-n}(\bar\xi_{V_N\ba V_n})\right)+\sum_{i=1}^n \P^{(1)}
\left(\frac{N-n}{N}L_{N-n}
(\bar\xi_{V_N\ba V_n}),\d_{\s_i}\right)\cr
+&o\left(\frac{1}{N}\right),
\end{split}
\end{equation}
where $o\left(\frac{1}{N}\right)$ is a result of \eqref{differentiability} of Definition \ref{mfinteraction} 
where we have taken $t=\frac{1}{N}$.
Finally, by putting this expression of $N\P\left(L_{N}(\bar\xi_{V_N})\right)$ \eqref{decom} into the expression
of\\ $\mu'_{n,N}(d\eta_{V_n}|\bar\eta_{V_N\ba V_n})$ \eqref{defin} and multiplying the resulting expression
by 
\begin{equation}
\dfrac{\prod_{j=n+1}^N \int \a(d\s_j) \exp\Bigl( -N\P\bigl(\pi_1\bar\nu_{N,n}\bigr)
	\Bigr) k(\s_j,\bar\eta_j)}{\prod_{j=n+1}^N \int \a(d\s_j)\exp\Bigl( -N\P\bigl(\pi_1\bar\nu_{N,n}\bigr)
	\Bigr) k(\s_j,\bar\eta_j)}.
\end{equation}
conclude the proof of the lemma.
\begin{flushright}
$\Cox$
\end{flushright}

We now state the infinite-volume ($N\rightarrow\infty$) version of Lemma \ref{mainlemma}.
 A sufficient condition for 
the existence of the finite-volume conditional distributions with infinite-volume 
$\eta$-conditioning is provided. This sufficient condition is  the unicity of the global 
minimizers of the function $J_{\nu'}$, $\nu'\in\PP(S')$.

\begin{thm} \label{maintheorem}
For each  $\nu'\in\PP(S')$, let $J_{\nu'}:M_{\nu'}\rightarrow \R\cup\{+\infty\}$ be as defined in 
\eqref{slayerrf}. Suppose further  that for a given $\nu'\in\PP(S')$ 
$J_{\nu'}$ has a unique global minimizer $\td\nu^*\in M_{\nu'}$.\\
I)   Then
\begin{equation}\begin{split}\label{gaman}
\g'_{n}(d\eta_{V_n}|\nu')=&\lim_{N\rightarrow \infty}\mu'_{n,N}(d\eta_{V_n}|\bar\eta_{V_N\ba V_n})=
\prod_{i\in V_n} \g'_{1}(d\eta_i|\nu'),\quad \text{with}\cr
 \g'_{1}(d\eta_i|\nu')=&
\dfrac{\displaystyle{ \int_{S}\a(d\s_i)e^{-\P^{(1)}\left(\pi_1\td\nu^*,\d_{\s_i}\right)}
k(\s_i,\eta_i)
\a'(d\eta_i)}}{\displaystyle{ \int_{S}\a(d\s_i) e^{-\P^{(1)}\left(\pi_1\td\nu^*,\d_{\s_i}\right)}
}}.
\end{split}
\end{equation}
II) If $J_{\l'}$ has a unique global minimizer $\td\l^*$ for all $\l'$ in a weak neighborhood of $\nu'$, then
 $\g'_{n}(d\eta_{V_n}|\nu')$ is weakly continuous at $\nu'$ as a function of the conditioning $\nu'\in\PP(S')$.
\end{thm}

\begin{remark}
\textbf{1. } The $\nu'$ dependence of the expression for $\g'_{1}(d\s'_1|\nu')$ is hidden 
in $\P^{(1)}$, via the probability measure $\pi_1\td\nu^*$.\\
\textbf{2. } Theorem \ref{maintheorem}  provides a sufficient condition for the transformed 
system to be Gibbs (in the sense of Definition \ref{gibbs}), namely the unicity of global 
minimizers of $J_\cdot$. Thus the problem of determining whether the transformed system is 
Gibbs or not is then translated into the corresponding problem of studying the global minimizers 
of $J_{\cdot}$. 
\end{remark}
\textbf{Proof of Theorem \ref{maintheorem}:} \\
I) The proof follows by way of the form of the finite-volume conditional distribution given in 
Lemma \ref{mainlemma} and the hypothesis that the function $J_{\nu'}$ has a unique global minimizer. Because the
leading term in the large $N$ asymptotic of $\mu'_{n,N}(d\eta_{V_n}|\bar\eta_{V_N\ba V_n})$ is governed
by the global minimizers of $J_{\nu'}$.\\\\
II) By hypothesis we get continuity for free by the continuity properties of the mean-field interaction 
$\P$.
\begin{flushright}
$\Cox$
\end{flushright}

\section{Gibbsianness of transformed systems and the contraction map theorem }\label{s4}
This section is devoted for studying the minimizers of the function $J_\cdot$ for some special class 
of initial interactions $\P$. Up to this point all topological considerations have been w.r.t. the 
weak topology, i.e. the weak topology is sufficient to study   Gibbs measures and Gibbs 
properties of transforms of Gibbs measures for mean-field models. We now consider another topology
on the spaces of measures  which is stronger than the weak topology.
This topology is the one induced by the total variational metric. Continuity in this new topology 
implies the continuity w.r.t. the weak topology. All topological considerations  for the interactions 
we consider in this section shall be w.r.t. the variational topology. Additionally, we also impose   further 
smoothness requirements  on the initial interactions $\P$ other than those given in Definition \ref{mfinteraction}. 
All these restrictions on the interactions are required to derive explicit continuity estimates on the CFLPK's and
consequently on the transformed kernels.
 
To be  precise, we  consider interactions $\P$ that are given by
\begin{equation}
\P(\nu)=F(\nu[g_1],\ldots,\nu[g_l]),
\end{equation}
where $g_i$ are some fixed bounded non-constant  real-valued measurable functions defined  on 
$S$, $l\geq1$ and $F:\R^l\rightarrow \R$ is some twice continuously differentiable 
function (e.g. if $F$ is a polynomial). In the following we will write $g=(g_1,
\cdots,g_l)$ and $\nu[g]=\bigl(\nu[g_1],\cdots,\nu[g_l]\bigr)$. By setting 
$m_j=m_j(\nu)=\nu[g_j]$, we have for this choice of interaction that
\begin{equation}
\P^{(1)}(\nu,\d_{\s_i})=\sum_{j=1}^l F_{j}(\nu[g_1],\ldots,\nu[g_l])\bigr)g_j(\s_i),
\end{equation}
where  $F_{j}(m)=\frac{\partial }{\partial m_j}F(m)$ and $m=m(\nu)=\nu[g]$. We also set 
$F_{ju}(m)=\frac{\partial^2}{\partial m_j\partial m_u}F(m)$.
Additionally, we assume that   $g$ is Lipschitz-function from $S$  to $\R^l$,
with  Lipschitz-norm 
\begin{equation}\begin{split}\label{Lipnorm}
 &\Vert g\Vert_{d,2} =\sup_{\s_i\neq \bar\s_i}\frac{\| g(\s_i)-g(\bar \s_i)\|_2}
 {d(\s_i,\bar\s_i)},
 \end{split}
\end{equation}
where $d$ is the metric on $S$. We also denote by $\d (g)$  the sum of the oscillations of the components of 
$g$, i.e.
\begin{equation}\label{deltag}
\d(g)=\sum_{j=1}^l \d(g_j).
\end{equation}
For any $g$ satisfying the above conditions we set
\begin{equation}\label{dbarg}
D_{g}=\overline{\left\{\nu[g]:\nu\in\PP(S)\right\}}.
\end{equation}
Note that $D_g$ is compact subset of $\R^l$ by the boundedness of $g$.
In the sequel we will  write $\|\partial^2F\|_{\text{max},\infty}$ for the supremum of the matrix max norm
of the Hessian $\partial^2 F$. i.e.
\begin{equation}\label{maxnorm}
\|\partial^2F\|_{\text{max},\infty}=\sup_{m\in D_g}\|\partial^2 F\big(m\big)\|_{{\rm max}},\quad 
\text{where}\quad \|\partial^2 F\big(m\big)\|_{{\rm max}}=\max_{1\leq i,j\leq l} \big| F_{ij}(m)\big|.
\end{equation}
Furthermore, we also set 
\begin{equation}
\hat{\d}\left(\P^{(1)}\right)=\sup_{m\in D_g}\sup_{\s_i,\bar\s_i\in S}\Big|\sum_{j=1}^l F_j(m)\Big(
g_j(\s_i)-g_j(\bar\s_i)\Big)\Big|.
\end{equation}
 Up to this point one may wonder whether the class of interactions we are considering in this section 
 has any physical relevance. Indeed, it contains important  mean-field interactions like 
 the Curie-Weiss interactions, liquid crystal interactions, sums of  ``$p$-spin'' interactions  
 ect., that have featured prominently in the literature.  
   
\subsection{Lipschitz Continuity of the CFLPK and Gibbsianness of Transformed System} 
We have already seen from the remark below equation \eqref{rckernel} that the CFLPK 
$\g_\cdot$ is weakly continuous. In this subsection we show, however that the CFLPK
is Lipschitz continuous w.r.t.  the variational metric (defined below). We write 
 \begin{equation}\begin{split}\label{vardist}
\|\nu-\bar \nu\|&=\sup_{|\phi|\leq1}|\nu(\phi)-\bar\nu(\phi)|\cr
&=\sup_{\bar\phi}\dfrac{\big|\nu(\bar\phi)-\bar\nu(\bar\phi)\big|}{\d
(\bar\phi)}\quad \text{where}\quad
\d(\bar\phi)=\sup_{\s_i\neq \bar\s_i}\big|\bar\phi(\s_i)-\bar\phi(\bar\s_i)\big|,
\end{split}
\end{equation} 
for the variational distance between the probability measures $\nu$ and $\bar\nu$ 
where the supremum are respectively taken  over all  measurable real-valued functions 
$\phi$ with 
$|\phi|\leq1$ and bounded non-constant measurable real-valued functions on $S$. 
The variational distance  can also be define by the following consideration:
The signed measure $\nu-\bar\nu$ has respectively  $(\nu-\bar\nu)^+$ and $(\nu-\bar\nu)^-$ as
the positive and negative parts of its  Jordan decomposition. But the fact that 
 $(\nu-\bar\nu)(S)=0$, implies that  $(\nu-
\bar\nu)^+(S)=(\nu-\bar\nu)^-(S)$ leading to the definition of the variational 
distance between $\nu$ and $\bar\nu$ as the one-half of the total variation of
$(\nu-\bar\nu)$, i.e.
 \begin{equation}
 \|\nu-\bar\nu\|=(\nu-\bar\nu)^+(S)=(\nu-\bar\nu)^-(S).
 \end{equation}

Before we state our  first result in  this section let us fix further notations.
We set
\begin{equation}\begin{split}
C(F,g)&=2\|\partial^2F\|_{\text{max},\infty}\; \d(g)\Vert g\Vert_{d,2}
\exp\Bigl(\frac{\hat{\d}\big(\P^{(1)}\big)}{2}\Bigr)\quad \text{and}\cr
\rho_{\a}(k)&=\sup_{\eta_i\in S'}\inf_{a_i\in S}\Bigl(\int_S d^2(\s_i,a_i)k(\s_i,\eta_i)
\a(d\s_i)\Bigr)^{\frac{1}{2}}.
\end{split}
\end{equation}

\begin{thm}\label{lipcont}
For any $\nu'\in\PP(S)$ and each pair $\td\nu_1,\td\nu_2\in M_{\nu'}$, the CFLPK 
satisfies
\begin{equation}\begin{split}\label{boundk}
&\|\g_{\nu'}(\cdot|\td\nu_1)-\g_{\nu'}(\cdot|\td\nu_2)\|\leq L \|\td\nu_1-\td\nu_2
\|\quad \text{where}\quad
L=L(F,g,k)= C(F,g)\rho_{\a}(k).
\end{split}
\end{equation}

\end{thm}
The above theorem says that for each $\nu'\in\PP(S')$, the CFLPK  $\g_{\nu'}$ is  Lipschitz continuous 
on $M_{\nu'}$ with Lipschitz constant $L$.
\begin{remark}\textbf{1. }
The quantity $\rho_\a(k)$ is the (metric-space version of) standard deviation 
of the single-site "posterior distribution" $K(d\s_i|\eta_i)$, when we take supremum 
over the possible observations $\eta$. So, it describes the worst-$\eta$ scenario 
of the typical size of fluctuations in the initial configurations which have led to $\eta$. 
The constant $L$ factorizes into two constants reflecting the idea of  "nature $C(F,g)$ 
versus nurture $\rho_\a(k)$".\\
\textbf{2. }  Set 
 \begin{equation}\label{varalpha}
 \rho_\a=\inf_{a_i\in S}\Bigl(\int \a(d\s_i)d^2(\s_i,a_i)\Bigr)^{\frac{1}{2}},
 \end{equation}
 i.e. $\rho_\a$ is the metric space version of the standard deviation of $\a$.
  Then the initial kernel $\g_1$ for the interactions considered in this section is also
  Lipschitz continuous, i.e. for any pair $\nu_1,\nu_2\in\PP(S)$ we have   
\begin{equation}\label{lipin}
\|\g_1(\cdot|\nu_1)-\g_1(\cdot|\nu_2)\|\leq \hat L \|\nu_1-\nu_2\|,\quad {where}\quad \hat L=
C(F,g)\rho_\a.
\end{equation}
$\hat L$ is the "Dobrushin's constant" for the initial mean-field model.  
\end{remark}
\textbf{Proof of Theorem \ref{lipcont}:}
Take a measurable map $f:S\times S'\rightarrow \R\;$  with $\;|f|\leq 1$. Also    for any pair 
$\td\nu_1,\td\nu_2\in M_{\nu'}$  and  any $0\leq s\leq1$ we define   $\td\nu_s:=s\td\nu_1+
(1-s)\td\nu_2$. Then we have 
\begin{equation}\begin{split}\label{prooflip}
\Big|\g_{\nu'}(f|\td\nu_1)-\g_{\nu'}(f|\td\nu_2)\Big|&=\Big|\int \nu'(d\eta_i)\int\a(d\s_i)
f(\s_i,\eta_i)\int_0^1 ds \frac{d}{ds}h_{\td\nu_s}(\s_i,\eta_i) \Big|,\;\text{where}\cr
h_{\td\nu_s}(\s_i,\eta_i)&=\dfrac{e^{-\P^{(1)}(\pi_1\td\nu_s,\d_{\s_i})}k(\s_i,\eta_i)}
{\int e^{-\P^{(1)}(\pi_1\td\nu_s,\d_{\s_i})}k(\s_i,\eta_i)\a(d\s_i)}.
\end{split}
\end{equation}
We also set $\l_s[\eta_i](d\s_i)=h_{\td\nu_s}(\s_i,\eta_i)\a(d\s_i)$.
Now using the form of the interactions considered in this section, it is not hard to deduce that
\begin{equation}
\frac{d}{ds}\P^{(1)}(\pi_1\td\nu_s,\d_{\s_i})=\sum_{j=1}^l\sum_{u=1}^l F_{ju}\big(m(\pi_1\td\nu_s)\big)
g_j(\s_i)\int\pi_1(\td\nu_1-\td\nu_2)(d\bar\s_i)g_u(\bar\s_i).
\end{equation}
 This and further computations yields the following expression for $\frac{d}{ds}h_{\td\nu_s}$:
\begin{equation}\begin{split}\label{direv}
\frac{d}{ds}h_{\td\nu_s}(\s_i,\eta_i)&=-\sum_{j=1}^l\sum_{u=1}^l F_{ju}\big(m(\pi_1\td\nu_s)\big)\psi_u
h_{\td\nu_s}(\s_i,\eta_i)
\Big(g_j(\s_i)-\l_s[\eta_i](g_j)\Big),\;\;\text{where}\cr
\psi_u&=\int \pi_1(\td\nu_1-\td\nu_2)(d\bar\s_i)g_u(\bar\s_i).
\end{split}
\end{equation}
Observe from  \eqref{vardist} that 
  \begin{equation}\begin{split}
 \big| \psi_u\big|&=\d(g_u)\dfrac{\Big|\int
 \nu'(d\eta_i)\int \big(\td\nu_1(d\s_i|\eta_i)-
 \td\nu_2(d\s_i|\eta_i)\big)g_u(\s_i)\Big|}{\d(g_u)}\cr
 &\leq \d(g_u)\|\td\nu_1-\td\nu_2\|.
  \end{split}
  \end{equation}
 Putting all these together we arrive at 
\begin{equation}\begin{split}\label{go20}
&\Big|\g_{\nu'}(f|\td\nu_1)-\g_{\nu'}(f|\td\nu_2)\Big|\leq\|\td\nu_1-\td\nu_2\| \sum_{j=1}^l
\sum_{u=1}^l
 \d(g_u)\int_0^1ds\big|F_{ju}\big(m(\pi_1\td\nu_s)\big)\big|\int\nu'(d\eta_i)\phi(s,\eta_i),\cr
 &\text{where}\cr
 &\phi(s,\eta_i)=
\int\a(d\s_i)h_{\td\nu_s}(\s_i,\eta_i)
\big|g_j(\s_i)-\l_s[\eta_i](g_j)\big|.
\end{split}
\end{equation} 
By adding and subtracting $g_j(a_i)$ (for any arbitrary $a_i\in S$) from the term $g_j(\s_i)-\l_s[\eta_i](g_j)$
 in the definition
of $\phi(s,\eta_i)$ and applying the triangle inequality we arrive at the following:
 \begin{equation}\begin{split}
 \phi(s,\eta_i)\leq
2\int\a(d\s_i)h_{\td\nu_s}(\s_i,\eta_i)
\big|g_j(\s_i)-g_j(a_i)\big|.
  \end{split}
  \end{equation}
Further, it follows from H\"older's inequality that 
\begin{equation}\begin{split}
 \int\a(d\s_i)h_{\td\nu_s}(\s_i,\eta_i)
\big|g_j(\s_i)-g_j(a_i)\big|
&\leq\left(\int\a(d\s_i)h_{\td\nu_s}(\s_i,\eta_i)
\big(g_j(\s_i)-g_j(a_i)\big)^2\right)^{\frac{1}{2}}.
\end{split}
\end{equation} 
Now by replacing $\;\big|F_{ju}\big(m(\pi_1\td\nu_s)\big)\big|\;$ with  $\;\|\partial^2
F\|_{\text{max},\infty}\;$ in \eqref{go20} and  using the fact the 
square root function is concave we obtain
\begin{equation}\begin{split}\label{go21}
\Big|\g_{\nu'}(f|\td\nu_1)-\g_{\nu'}(f|\td\nu_2)\Big|&\leq 2\|\td\nu_1-\td\nu_2\|\;\d(g)
\|\partial^2F\|_{{\rm max},\infty} 
\int_0^1ds\int\nu'(d\eta_i)\chi(s,\eta_i,a_i)\quad
 \text{where}\cr
 \chi(s,\eta_i,a_i)&=\Bigl(\int\a(d\s_i)h_{\td\nu_s}(\s_i,\eta_i)\sum_{j=1}^l\big(g_j(\s_i)-
g_j(a_i)\big)^2\Bigr)^{\frac{1}{2}}.
\end{split}
\end{equation}
Note further that the $a_i$ appearing in $\chi$ is chosen independent of all the parameters 
in the model, so taking the infimum over $a_i$ will have no influence on our estimates. In view of 
this observation, replacing $\chi(s,\eta_i,a_i)$ with $\inf_{a\in S}\chi(s,\eta_i,a_i)$ will have
no effect on the inequality in \eqref{go21}. Furthermore, it follows from the Lipschitz property
of $g$ and the fact that 
$h_{\td\nu_s}(\s_i,\eta_i)\leq e^{\hat{\d}(\P^{(1)})}k(\s_i,\eta_i)$ that  
\begin{equation}\begin{split}\label{go22}
\inf_{a_i\in S}\chi(s,\eta_i,a_i)&\leq \|g\|_{d,2}e^{\frac{\hat{\d}(\P^{(1)})}{2}}
\sup_{\eta_i\in S'} \inf_{a_i\in S} \Bigl(\int\a(d\s_i)k(\s_i,\eta_i)d^2(\s_i,a_i)\Bigr)^{\frac{1}{2}}.
\end{split}
\end{equation} 
Putting the bound in \eqref{go22} into the  bound in \eqref{go21} yields the desired result.
\begin{flushright}
$\Cox$
\end{flushright}
  
 Observe from Theorem \ref{lipcont} that if the constant $L<1$ then $\g_{\nu'}$ defines a contraction 
 map from $M_{\nu'}$ to itself. This is because the variational distance turns the set $M_{\nu'}=
 \pi_2^{-1}(\{\nu'\})$ into a complete metric space by continuity of the map $\pi_2$ under the 
 variational topology. Thus the CFLPK admits a unique  consistent probability measure and 
 consequently the existence of a unique global minimizer for 
 $J_\cdot$, since the minimizers of $J_\cdot$ are contained in the set of consistent probability 
 measures for the CFLPK. This then implies that the transformed system  is Gibbs. 
 
 The next item  on our list of tasks is the investigation of how (in the regime $L<1$) the unique 
 consistent probability measure 
 $\td\nu^*$ for  $\g_{\nu'}$ behaves w.r.t. $\nu'\in\PP(S')$.
 Indeed we show in the proposition below that $\td\nu^*$ depends continuously on $\nu'$. 
 
 \begin{prop}\label{conspro}
 Suppose the constant $L<1$, then under the variational metric the unique consistent 
 probability measure $\td\nu^*$ for 
 $\g_{\nu'}$ is  Lipschitz continuous w.r.t. $\nu'$ and has the Lipschitz norm $L_1=4\;L$.
 \end{prop}
 
 \begin{remark}
 The constant $L$ here is comparable with the uniform bound on the Dobrushin constant $c'[\eta]$ for 
 the restricted constrained first layer model considered in \cite{KULOP}. Due to this, it is of 
 interest to obtain  the Lipschitz constant with the factor $\frac{1}{1-L}$ reminiscent to
 the upper bound on the row sums of the Dobrushin matrix $D$. Thus for $L<1$ we can also have  
 the constant $\bar L_1$ given by
 \begin{equation}\label{secondconst}
 \bar L_1= \dfrac{L_1}{1-L}, \quad \text{with} \quad L_1<\bar L_1.
 \end{equation}
 \end{remark}
 \textbf{Proof of Proposition \ref{conspro}: } 
 Let $\nu'_1,\nu'_2\in\PP(S')$, then the assertion of the proposition follows by showing that,
 \begin{equation}
 \|\g_{\nu'_1}-\g_{\nu'_2}\|\leq L_1\|\nu'_1-\nu'_2\|,
 \end{equation}
 since if $\td\nu^*_i$ is the unique consistent probability measure for $\g_{\nu'_i}$, 
 then $\td\nu^*_i=\g_{\nu'_i}(\cdot|\td\nu^*_i)$ for $i=1,2$. Observe for any
 measurable function on $S\times S'$ with $|f|\leq 1$ that
 \begin{equation}\begin{split}
 \Big|\g_{\nu'_1}(f|\td\nu^*_1)-\g_{\nu'_2}(f|\td\nu^*_2)\Big|&=\Big|\int (\nu'_1-\nu'_2)(d\eta_i)
 \a(d\s_i)f(\s_i,\eta_i)\Bigl(h_{\td\nu^*_1}(\s_i,\eta_i)-h_{\td\nu^*_2}(\s_i,\eta_i)\Bigr)\Big|\cr
 &\leq \Big|\int (\nu'_1-\nu'_2)^+(d\eta_i)
 \a(d\s_i)f(\s_i,\eta_i)\Bigl(h_{\td\nu^*_1}(\s_i,\eta_i)-h_{\td\nu^*_2}(\s_i,\eta_i)\Bigr)\Big|\cr
 &+\Big|\int (\nu'_1-\nu'_2)^-(d\eta_i)
 \a(d\s_i)f(\s_i,\eta_i)\Bigl(h_{\td\nu^*_1}(\s_i,\eta_i)-h_{\td\nu^*_2}(\s_i,\eta_i)\Bigr)\Big|,
 \end{split}
 \end{equation}
 where $\;(\nu'_1-\nu'_2)^+\;$ and  $\;(\nu'_1-\nu'_2)^-\;$ are respectively the positive and the
 negative parts of the Jordan decomposition of the signed measure $\nu'_1-\nu'_2$ and $h_{\td\nu}$
 is as given in \eqref{prooflip}. It follows from 
 the definition of the variational distance between two probability measures that   
  
 \begin{equation}\begin{split}\label{58}
 \Big|\g_{\nu'_1}(f|\td\nu^*_1)-\g_{\nu'_2}(f|\td\nu^*_2)\Big|
 &\leq 2\|\nu'_1-\nu'_2\|\sup_{\eta_i\in S'}\int\a(d\s_i)\Big|h_{\td\nu^*_1}(\s_i,\eta_i)-
 h_{\td\nu^*_2}(\s_i,\eta_i)\Big|,
 \end{split}
 \end{equation} 
 since we have chosen $f$ to be  such that $\;|f|\leq 1$. Now we proceed by setting \\$\;\td\nu^*_s=
 s\td\nu^*_1+(1-s)\td\nu^*_2\;$ with $\;0\leq s\leq 1\;$ and
 observing that
 \begin{equation}\begin{split}
\int\a(d\s_i)\Big|h_{\td\nu^*_1}(\s_i,\eta_i)-h_{\td\nu^*_2}(\s_i,\eta_i)\Big|&=
\int\a(d\s_i)\Big|\int_0^1 ds\frac{d}{ds}h_{\td\nu^*_s}(\s_i,\eta_i)\Big|\cr
&\leq \int_0^1 ds \int \a(d\s_i)\Big|\frac{d}{ds}h_{\td\nu^*_s}(\s_i,\eta_i)\Big|.
 \end{split}
 \end{equation} 
 But we know from the proof of Theorem \ref{lipcont} that 
 \begin{equation}\begin{split}\label{60}
 \int_0^1 ds \int \a(d\s_i)\Big|\frac{d}{ds}h_{\td\nu^*_s}(\s_i,\eta_i)\Big|\leq 
 C(F,g)\|\td\nu^*_1-\td\nu^*_2\|\inf_{a_i\in S}\Bigl(\int\a(d\s_i)k(\s_i,\eta_i)d^2(\s_i,a_i)
 \Bigr)^{\frac{1}{2}}.
 \end{split}
 \end{equation} 
 Now by replacing $\|\td\nu^*_1-\td\nu^*_2\|$ by $2$  and substituting the bound in \eqref{60} 
 into the inequality in \eqref{58} conclude the proof of the proposition.
\begin{flushright}
$\Cox$
\end{flushright}
 Having disposed of the continuity estimates for the CFLPK, we now turn  our attention to study  
 the corresponding continuity estimates for the single-site kernel $\g'_1$ for the transformed 
 system  in the regime where $L<1$. Since in this regime the CFLPK has a unique consistent 
 probability measure and consequently providing a sufficient condition for  the function $J_\cdot$ 
 to have a unique global minimizer as is required by Theorem \ref{maintheorem}. 
%
 \subsection{Continuity estimates for $\g'_1$}
 A sufficient condition given in Theorem \ref{maintheorem} for the existence and continuity of the 
 finite-volume kernels $\g'_n$ for the transformed system is the unicity of the global minimizer of 
 the function $J_{\nu'}$. This sufficient condition holds  if the CFLPK has a unique
 consistent probability. Our main result in this 
 subsection is the following theorem concerning the continuity estimate for $\g'_1$:
 
 \begin{thm}\label{contestimates}
 Suppose that $L<1$, then
 under the variational metric on $\PP(S')$ the single-site kernel $\g'_1(\cdot|\nu')$
 is Lipschitz continuous w.r.t. $\nu'$  with Lipschitz constant $L_2$ given by
 \begin{equation}\label{L2}
 L_2= L_1 \hat L.
 \end{equation} 
 \end{thm}

\begin{remark}
\textbf{1. }  The Lipschitz constant $L_2$ for the transformed system factorizes into the product
of the Lipschitz constants $L$ and $\hat L$ respectively for the CFLPK $\g_\cdot$ and the initial kernel 
 $\g$.\\
 \textbf{2. } Observe from the remark below Proposition \ref{conspro} that we can also have the Lipschitz constant
 $\bar L_2$ given by 
 \begin{equation}
 \bar L_2=\bar L_1 \hat L=\dfrac{4L\hat L}{1-L} \quad \text{with}\quad L_2<\bar L_2.
 \end{equation}
 \end{remark}
 \textbf{Proof of Theorem \ref{contestimates}: }
 As usual let us take $\td\nu^*_1$ and $\td\nu^*_2$ as the  unique consistent probability
 measures for the CFLPK corresponding to  $\nu'_1$ and $\nu'_2$ respectively. Again set
 $\td\nu^*_s=s\td\nu^*_1+(1-s)\td\nu^*_2$ for $0\leq s\leq 1$. It follows from \eqref{gaman} after 
 taking a measurable map $f:S'\rightarrow\R$ with $|f|\leq1$ and setting 
 \begin{equation}
 \hat h_{\td\nu^*_s}(\s_i,\eta_i)=\dfrac{k(\s_i,\eta_i)e^{-\P^{(1)}(\pi_1\td\nu^*_s,\d_{\s_i})}}
 {\int e^{-\P^{(1)}(\pi_1\td\nu^*_s,\d_{\hat\s_i})}\a(d\hat\s_i)}
 \end{equation}
that
\begin{equation}\begin{split}
&\g'_1(f|\nu'_1)-\g'_1(f|\nu'_2)=\int \a'(d\eta_i)f(\eta_i)\int\a(d\s_i)\int_0^1ds\frac{d}{ds}\hat
 h_{\td\nu^*_s}(\s_i,\eta_i)\quad\text{with}\cr
 &\frac{d}{ds}\hat h_{\td\nu^*_s}(\s_i,\eta_i)=\cr
 &-\sum_{j=1}^l\sum_{u=1}^l F_{ju}\big(m(\pi_1\td\nu^*_s)\big)
 \int \pi_1(\td\nu^*_1-\td\nu^*_2)(d\bar\s_i)g_u(\bar\s_i)\hat h_{\td\nu^*_s}(\s_i,\eta_i)\Bigl(g_j(\s_i)
 -\vartheta_s(g_j)\Bigr)\cr
 &\text{and where}\quad \vartheta_s(g_j)=\int \a'(d\eta_i)\a(d\s_i)\hat h_{\td\nu^*_s}(\s_i,\eta_i)g_j(\s_i).
 \end{split}
 \end{equation} 
 Therefore it follows from our previous considerations that 
 \begin{equation}\begin{split}\label{66}
\Big|\g'_1(f|\nu'_1)-\g'_1(f|\nu'_2)\Big|&\leq\|\partial^2F\|_{{\rm max},\infty}\;\d(g)
\|\td\nu^*_1-\td\nu^*_2\|\times\cr
&\sum_{j=1}^l\int_0^1 ds \int \a'(d\eta_i)\a(d\s_i)\hat
 h_{\td\nu^*_s}(\s_i,\eta_i)\Bigl|g_j(\s_i)
 -\vartheta_s(g_j)\Bigr|.
 \end{split}
 \end{equation} 
 By adding and subtracting $g_j(a_i)$ from the term $g_j(\s_i)-\vartheta_s(g_j)$ and applying 
the  triangle inequality 
 we  obtain
\begin{equation}\begin{split}\label{69}
\int \a'(d\eta_i)\a(d\s_i)\hat h_{\td\nu^*_s}(\s_i,\eta_i)\Bigl|g_j(\s_i)-\vartheta_s(g_j)\Bigr|&\leq
 2\int\a'(d\eta_i)\a(d\s_i)\hat h_{\td\nu^*_s}(\s_i,\eta_i)\Bigl|g_j(\s_i)-g_j(a_i)\Bigr|
 \end{split}
 \end{equation}

Now it follows from H\"older's inequality and the facts that (1) the square 
root function is
concave (2) $\;\hat h_{\td\nu^*_s}(\s_i,\eta_i)\leq e^{\hat{\d}(\P^{(1)})}k(\s_i,\eta_i)\;$
and (3) the Lipschitz property of $g_j$  that 
 
 \begin{equation}\begin{split}\label{67}
 &\sum_{j=1}^l\int_0^1 ds \int \a'(d\eta_i)\a(d\s_i)\hat h_{\td\nu^*_s}(\s_i,\eta_i)\Bigl|g_j(\s_i)
 -\vartheta_s(g_j)\Bigr|\cr
 &\leq 2\int_0^1 ds \inf_{a_i\in S}\sum_{j=1}^l\Bigl(\int \a'(d\eta_i)\a(d\s_i)\hat
 h_{\td\nu^*_s}(\s_i,\eta_i)\bigl(g_j(\s_i)-g_j(a_i)\bigr)^2\Bigr)^{\frac{1}{2}}\cr
 &\leq 2\int_0^1 ds \inf_{a_i\in S}\Bigl(\int \a'(d\eta_i)\a(d\s_i)\hat
 h_{\td\nu^*_s}(\s_i,\eta_i)\sum_{j=1}^l\bigl(g_j(\s_i)-g_j(a_i)\bigr)^2\Bigr)^{\frac{1}{2}}\cr
 &\leq 2\|g\|_{d,2}e^{\frac{\hat{\d}(\P^{(1)})}{2}}\inf_{a_i\in S}\Bigl(\int \a(d\s_i)d^2(\s_i,a_i)\Bigr)^{\frac{1}{2}}.
 \end{split}
 \end{equation} 
 Finally by putting this bound in \eqref{67} into \eqref{66} and noting from Proposition \ref{conspro}
 that $\|\td\nu^*_1-\td\nu^*_2\|\leq L_1 \|\nu'_1-\nu'_2\|$ follows the proof. 
\begin{flushright}
$\Cox$
\end{flushright}

\section{ Examples}\label{s5}
We now present two examples  for the class of models discussed in the above section.  In the first 
example we consider  specific forms of the functions $F$ and $g$, and a specific form of the joint 
a priori measure $K$. The second example is about   general forms of $F$ and $g$, and a 
specific form of the joint a priori measure $K$.

\subsection{Short time Gibbsianness of rotator mean-field models  under diffusive time evolution}
 The first example we consider is the Curie-Weiss rotator  model 
under site-wise independent diffusive time-evolution. Here the single-site spin space for both the initial and
the transformed systems are the same, i.e. 
$S=S'=S^{q-1}$, where $S^{q-1}$ is  the sphere in the q-dimensional Euclidean space with $q\geq 2$. The 
interaction for the initial system is given by 
\begin{equation}\label{cweiss}
\P(\nu)=F\bigl(\nu[\s^1_i],\cdots,\nu[\s^q_i]\bigr)=-\dfrac{\beta \sum_{j=1}^q \nu[\s^j_i]^2}{2},
\end{equation}
where $g_{j}(\s_i)=\s_i^j$ is the $j$th coordinate of the point $\s_i\in S^{q-1}$ and 
 $l=q$. 

Next let $K$  be given by 
 \begin{equation}
 K(d\s_i,d\eta_i)=K_t(d\s_i,d\eta_i)=k_t(\s_i,\eta_i)\a_0(d\s_i)\a_0(d\eta_i),
  \end{equation}
where $\a_0$ is the equidistribution on  $S^{q-1}$ and
$k_t$ is the heat kernel on the sphere, i.e.
\begin{eqnarray}
\Big(e^{\D t}\varphi\Big)(\eta_i)=\int \a_o(d\s_i)k_t(\s_i,\eta_i)\varphi(\s_i),
\end{eqnarray} 
where $\D$ is the Laplace-Beltrami operator on the sphere and $\varphi$ is any test function.  
$k_t$ is also called the {\em  Gauss-Weierstrass kernel}. 
For more background on the heat-kernel on Riemannian manifolds, see the introduction 
of \cite{ACDH04}. 

\begin{lem}\label{time-evol}
Assuming the above set-up the  Lipschitz constant $L$ for the CFLPK  is given by
\begin{equation}
L=L(q,\beta,t)=4\sqrt{2}q\beta e^\beta\Bigl(1-e^{-(q-1)t}\Bigr)^{\frac{1}{2}}.
\end{equation}
\end{lem}
\textbf{Proof Lemma \ref{time-evol}:}
The proof follows by just evaluating the terms appearing in $C(F,g)$ and $\rho_\a(k)$.
First observe that for the interaction considered in this subsection $\|\partial^2F\|_{\text{max},\infty}=\beta$,
$\;\d(g)=2q$, $\;\|g\|_{d,2}=1\;$ since $\;d(\s_i,\bar\s_i)=\|g(\s_i)-g(\bar\s_i)\|_2\;$ and 
$\;\hat{\d}(\P^{(1)})=2\beta$. So these give rise to
\begin{equation}
C(F,g)=C(q,\beta)=4 q\beta e^{\beta}.
\end{equation}
To obtain $\rho_\a(k)$ we argue as follows:
Denote by $u(t)$ the $q$-th component of a diffusion $\s_i(t)$ on the sphere
started at $u(0)=1$ (in the "north-pole") and  $\E$ the expectation w.r.t. the
corresponding diffusion. By reversibility we now  choose $\;a_i=\eta_i\;$ such that $\;a_i\;$ is the north pole . Then
it is not hard to see that 
$$d^2(\s_i(t),a_i)=2\big(1-u(t)\big).$$
This gives rise to
\begin{equation}\begin{split}\label{ehoch}
 &  \int \a_0(d\s) k_t(\s_i(t),a_i) d^2(\s_i(t),a_i) = 2 (1- \E u(t) ).
\end{split}
\end{equation}
It follows from the above and the rotation invariance of the diffusion on the sphere that to compute 
$\rho_\a(k)$ we only need 
the $q$th component of the diffusion which according to the Laplace-Beltrami operator on the sphere has generator of the 
form 
\begin{equation}
(1-u^2)\Bigl(\frac{d}{du}\Bigr)^2 -(q-1)u \frac{d}{du}.
\end{equation}
 This  generates the equation 
 $\frac{d}{dt}\E u_t = -(q-1)\E u_t $. Solving this equation  with the initial 
 condition $u(t=0)=1$ yields 
 \begin{equation}
\E u(t)=e^{-(q-1)t},
\end{equation}
 which leads to the desired expression for $\rho_\a(k)$, i.e.
\begin{eqnarray}
\rho_\a(k)&=&\sqrt{2}\Bigr(1-e^{-(q-1)t} \Bigl)^{\frac{1}{2}}\;\text{and\; therefore}\\
L&=&4\sqrt{2} q \beta\;e^{\beta}\Bigr(1-e^{-(q-1)t} \Bigl)^{\frac{1}{2}}\nonumber
\end{eqnarray}  
\begin{flushright}
$\Cox$
\end{flushright}

\begin{remark}
\begin{enumerate}
	\item The smallness of $L$ for this example emanates from  at least  two sources, namely; 
	small values of $t$ and $\beta$. That is if $\beta$ is small enough the system will be Gibbs at all times.
	But if we start with large $\beta$ then we hope to preserve Gibbsianness at only small values of $t$.
\item We also have for any arbitrary chosen $a_i$ that
\begin{equation}
\int_S d^2(\s_i,a_i) \a(d\s_i)=\int 2(1-\s\cdot a)\a_o(d\s)=2.
\end{equation} 
Thus we have for this example $\rho_\a=\sqrt{2}.$
\end{enumerate}
 \end{remark}

\subsection{ Local approximation and  preservation of Gibbsianness}
As our second example we start from an initial compact Polish space $S$ endowed
with a metric $d$ and an a priori measure $\a$. 
We consider general $F$ and $g$ defining the initial Hamiltonian. 

We partition the initial space $S$ into finitely or countably infinitely many disjoint 
Borel sets with  non-zero $\a$ measure indexed by the elements in $S'$, i.e. 
\begin{equation}
S=\bigcup_{\eta_i\in S'} S_{\eta_i},\quad \text{with} \quad \a(S_{\eta_i})>0\quad \text{for all}\quad \eta_i\in S'.
\end{equation}

We then consider the deterministic map $T:S\rightarrow S'$,
such that $T(\s_i)=\eta_i$ for all  $\s_i\in S_{\eta_i}$. That is, every point is mapped to the label 
of the class it belongs to. 
If we start with a finite initial space, this transformation is the so-called fuzzy-map 
which, when starting from an initial Potts model, was studied in \cite{HaggKu04}.  
In the present generality this example was studied in \cite{KULOP}, and we want to see here 
what the mean-field estimates of the present paper provide. 
Let us  formulate  the form of the Lipschitz constant $L$ for the CFLPK 
resulting from the local approximations  in the following lemma. 

\begin{lem}
Assume the set-up above, then the Lipschitz constant $L$ is given by
\begin{equation}
L=L(F,g,T)=C(F,g)\sup_{\eta_i\in S'}\a(S_{\eta_i})^{-\frac{1}{2}}\inf_{a_i\in S_{\eta_i}}\Bigl(\int \a_{|_{S_{\eta_i}}}
(d\s_i) d^2(\s_i,a_i)\Bigr)^{\frac{1}{2}}.
\end{equation}
\end{lem}
The proof of the above Lemma follows straight away from the definition of the constant $L$ and observing that
$\a(d\s_i)k(\s_i,\eta_i)=\frac{\a_{|_{S_{\eta_i}}}(d\s_i)}{\a(S_{\eta_i})}$.
 
Once again, the constant $L$ will be small either if the initial interaction is weak enough or the local 
approximation is fine enough.  For $L<1$, by the general Theorem \ref{contestimates},
this implies Gibbsianness and continuity estimates 
of the form \eqref{L2}.

 \section{Appendix}
 \subsection{Proof of Theorem \ref{minimizers}}
 The assertion that $J_{\nu'}$ attains its infimum on $M_{\nu'}$ trivially follows from the lower 
 semi-continuity of  $J_{\nu'}$ and compactness of $M_{\nu'}$. Now to proceed with the rest of the 
 proof  we take any minimizer $\td\nu^*$ of $J_{\nu'}$ which has the representation
 $\td\nu^*(d\xi_i)=\nu'(d\eta_i)\a(d\s_i)f_{\nu'}(\s_i|\eta_i)$. Then it remains to show that 
 $f_{\nu'}$  (1) is $\a\otimes\nu'$-a.s. strictly  positive and (2) takes the form \eqref{cdf} 
 $\a\otimes\nu'$-a.s..\\\\ 
(1)$\;$ We now proceed to  show the almost sure strict
positivity of the minimizing conditional $\a-$density $f_{\nu'}$. 
For each $\eta_i\in S'$ we set  $A_{\eta_i}=\left\lbrace \s_i\in S: f_{\nu'}(\s_i\vert\eta_i)=0 
\right\rbrace$ and denote $B$ by the set of $\eta_i$'s for which $\alpha(A_{\eta_i})>0$. 
The proof consists of establishing a contradiction that $\td\nu^*$ is not a minimizer of $J_{\nu'}$ 
whenever $B$ has a positive $ \nu'$ measure.  This is an adaptation of arguments found in \cite{LEW90}
 and references therein  modified to suit our case.

To be precise  let 
$b:S\rightarrow\R$ be a strictly positive measurable map with $b\leq1$ and for each $\eta_i\in B$ we define
a bounded measurable map $g(\cdot|\eta_i):S\rightarrow\R$ by
\begin{eqnarray}
g(\s_i\vert\eta_i)=\frac{\left(1_{A_{\eta_i}^{c}}f_{\nu'}(\cdot\vert\eta_i)\right)(\s_i)+\left(1_{A_{\eta_i}}b\right)
(\s_i)}{1+\int_{A_{\eta_i}} b(\hat\s_i)\a(d\hat\s_i)}.
\end{eqnarray} 
Now set $u_{\eta_i}=1+\int_{A_{\eta_i}} b(\hat\s_i) \a(d\hat\s_i)$  and further define for each $\eta_i\in S'$

\begin{equation}\begin{split}
p_{\e}(\cdot,\eta_i)&=\left\{
\begin{array}{rl}
\e g(\cdot\vert\eta_i)+(1-\e)f_{\nu'}(\cdot\vert\eta_i)=\frac{\e}{u_{\eta_i}}1_{A_{\eta_i}}b+
C_{\e,\eta_i}1_{A_{\eta_i}^{c}}f_{\nu'}(\cdot\vert\eta_i) & \text{if } \eta_i\in B\\\\
f_{\nu'}(\cdot\vert\eta_i) & \text{if } \eta_i\in B^c
\end{array} \right.
,\;\text{where}\cr
C_{\e,\eta_i}&=\e(\frac{1}{u_{\eta_i}}-1)+1
\end{split}
\end{equation} 
and  $\e\in[0,1]$. It is easy to check that  $1\leq u_{\eta_i}\leq 2$, which implies that
 $C_{\e,\eta_i}$ as well as $\log C_{\e,\eta_i}$ is uniformly bounded. Let us set 
$\td\nu^*_\e(d\xi_i)=\nu'(d\eta_i)\a(\s_i)p_{\e}(\s_i,\eta_i)$ and observe 
from the above that we can write the relative entropy of 
$\td\nu^*_\e$ w.r.t. $\a\otimes \nu'$ as; 
\begin{equation}\begin{split}
S(\td\nu^*_\e|\a\otimes \nu' )&=  \int_{B^c}\nu'(d\eta_i)S\left(\td\nu^*(\cdot|\eta_i)\big|\a\right) +
\int_{B}\nu'(d\eta_i)\bigg\{\frac{\e}{u_{\eta_i}}\int_{A_{\eta_i}}\a(d\s_i)b(\s_i)\log b(\s_i)\cr
&+\frac{\e}{u_{\eta_i}}\log\frac{\e}{u_{\eta_i}}\int_{A_{\eta_i}}\a(d\s_i)b(\s_i)+ C_{\e,\eta_i}
\int_{A^c_{\eta_i}}\a(d\s_i)f_{\nu'}(\s_i|\eta_i)\log f_{\nu'}(\s_i|\eta_i)\cr
&+C_{\e,\eta_i}\log C_{\e,\eta_i}
\int_{A^c_{\eta_i}}\a(d\s_i)f_{\nu'}(\s_i|\eta_i)\bigg\}.
\end{split}
\end{equation}
Now we define a function $h:[0,1]\rightarrow\R$ by
\begin{equation}\begin{split}
h(\e)=J_{\nu'}(\td\nu_\e^*)=S(\td\nu^*_\e|\a\otimes \nu' )+\P(\pi_1\td\nu^*_\e)-\td\nu^*_\e\left[\log k\right].
\end{split}
\end{equation}
It follows from the uniform boundedness of $C_{\e,\eta_i}$ and $\log C_{\e,\eta_i}$ , Lebesgue's dominated 
convergence theorem  and the continuity property of $\P^{(1)}$ that $h$ is continuously differentiable on $(0,1)$. 
Observe that $h(0)=J_{\nu'}(\td\nu^*)$ and one would expect $h$ to be decreasing as $\e\downarrow 0$, i.e. 
 $h(\e)-h(0)>0$ for $\e$ close to zero. But we will show that the converse of the above holds  if $B$ has 
positive $\nu'$ measure.
More precisely, differentiating $h$ we obtain for $\e\in(0,1)$
\begin{equation}\begin{split}
h'(\e)&= \int_B\nu'(d\eta_i) \bigg\{\frac{1}{u_{\eta_i}}\log \frac{\e}{u_{\eta_i}}\int\a(d\s_i)b(\s_i)+
\left(\frac{1}{u_{\eta_i}}-1\right)\log C_{\e,\eta_i}\int_{A^c_{\eta_i}}\a(d\s_i)f_{\nu'}(\s_i|\eta_i)\bigg\}\cr
&+\int_{B}\nu'(d\eta_i)\int_{A_{\eta_i}}\a(d\s_i)\left(g(\s_i|\eta_i)-f_{\nu'}(\s_i|\eta_i)\right)\P^{(1)}\left(
\pi_1\td\nu^*_\e,\d_{\s_i}\right) +C(B,b,f_{\nu'}),
\end{split}
\end{equation}
where $C(B,b,f_{\nu'})$ is a constant which depends on $B,\; b$ and $f_{\nu'}$ but independent of $\e$. 
Assuming 
$\nu'(B)>0$ then the limit 
$$\lim_{\e\downarrow 0}h'(\e)=-\infty$$
since the term 
\begin{equation}\begin{split}
\int_B\nu'(d\eta_i)\frac{1}{u_{\eta_i}}\log \frac{\e}{u_{\eta_i}}\int_{A_{\eta_i}}\a(d\s_i)b(\s_i)
\end{split}
\end{equation}
goes to negative infinity whiles the rest remains bounded.
This implies that
\begin{equation}\begin{split}
\lim_{\e\downarrow 0}\frac{h(\e)-h(0)}{\e}=-\infty
\end{split}
\end{equation}
giving rise to a contradiction since $\td\nu^*$ is a minimizer, hence the assumption than $\nu'(B)>0$
 is false. This concludes the proof of almost sure strict positivity of $f_{\nu'}$.\\\\
(2) Next we prove that the conditional $\a$-density of  any minimizer of $J_{\nu'}$ must satisfy the constrained 
mean-field equation \eqref{cdf}. 
Let $f_{\nu'}(\cdot\vert\cdot)$ be a conditional $\a$-density of a  minimizer $\td\nu^*$ of $J_{\nu'}$ and  set 
\begin{equation}\begin{split}\label{yeta}
Y_{\eta_i}&=\P^{(1)}(\pi_1\td\nu^*,\d_{\s_i})-\log k(\s_i,\eta_i)+\log f_{\nu'}(\s_i\vert\eta_i)\cr
&+\log\int\exp\left(-\P^{(1)}(\pi_1\td\nu^*,\d_{\hat\s_i})\right)k(\hat\s_i,\eta_i)\a(\hat\s_i).
\end{split}
\end{equation} 
We are now left  to show that $\a(Y_{\eta_i}=0)=1$ for $\nu'$ almost all $\eta_i\in S'$. 

The idea of the proof is again to assume the contrary, and arrive 
at the contradiction that  
a suitable perturbation of the conditional $\a$-density  $f_{\nu'}$ would have a lower value of $J_{\nu'}$.  
That is we assume  $\a(Y_{\eta_i}\neq 0)>0$ on some subset  $B$ of $S'$ with positive 
$\nu'$ measure. In a first step  this implies that for any $\eta_i\in B$ 
both inequalities $\a(Y_{\eta_i}\geq  0)> 0$ and $\a(Y_{\eta_i}<0)>0$ must be the case. 
Indeed, the assumption that  e.g. the second inequality is not true 
leads to a contradiction. To see this define for each $\eta_i\in B$  and $\delta>0$ 
\begin{equation}\begin{split}
A_{\eta_i}^{\delta}&=\Bigl\{ \s_i\in S: \psi_{\nu'}(\s_i|\eta_i) 
\geq\delta \Bigr\},\quad {with}\cr
\psi_{\nu'}(\s_i|\eta_i)&=f_{\nu'}(\s_i\vert\eta_i)-\frac{\exp\left(-\P^{(1)}
(\pi_1\td\nu^*,\d_{\s_i})\right)k(\s_i,\eta_i)}{\int\exp
\left(-\P^{(1)}(\pi_1\td\nu^*,\d_{\hat\s_i})\right)k(\hat\s_i,\eta_i)\a(\hat\s_i)}
\end{split}
\end{equation}
so that we would have $\lim_{\delta\downarrow 0}\a(A_{\eta_i}^{\delta})=1$.

Taking the  $\a$ integral of $\psi_{\nu'}(\s_i|\eta_i) $  yields 
then the contradiction 
\begin{equation}\begin{split}
0=&\int  \a(d\s_i)\psi_{\nu'}(\s_i|\eta_i)\geq 
  \int_{A_{\eta_i}^{\delta}}\a(d\s_i)\psi_{\nu'}(\s_i|\eta_i)\geq
 \delta\a(A_{\eta_i}^{\delta})>0
\end{split}
\end{equation}  for $\delta$ sufficiently  small. This and the $\a\otimes \nu'$-a.s. strict
positivity of the minimizing conditional $\a$-density $f_{\nu'}$ imply that $\a\big(1_{\{ Y_{\eta_i}
< 0\}}f_{\nu'}(\cdot\vert\eta_i)\big)>0$ and 
$\a\big(1_{\{ Y_{\eta_i}\geq 0\}}f_{\nu'}(\cdot\vert\eta_i)\big)>0$. Now we   set for each $\eta_i\in S'$
\begin{equation}\label{ceta}
C_{\eta_i}=\frac{\a\big(1_{\{ Y_{\eta_i}\geq 0\}}f_{\nu'}(\cdot\vert\eta_i)\big)}{\a\big(1_{\{ Y_{\eta_i}
< 0\}}f_{\nu'}(\cdot\vert\eta_i)\big)}
\end{equation} and define for each positive integer $n\in\N$ the set
  \begin{equation}
 B_{n}:=
 \Bigl\{\eta_i: C_{\eta_i}\in (n-1,n]\Bigr\}. 
 \end{equation}
 Observe that  the  $(B_n)_{n\geq1}$ is a partition of the set $B\subset S'$, i.e. $B=\bigcup_{n\geq1}B_{n}$ 
 and the $B_n$'s are pair-wise disjoint Borel subsets of $S'$.  We now consider a perturbation of the minimizing 
 conditional $\a$-density $f_{\nu'}(\cdot|\eta_i)$ whose form will be dependent on the choice of the transformed
 configuration $\eta_i$. More precisely, for each $\e\in[0,1]$ we consider the perturbed conditional 
 $\a$-density $p_\e(\s_i|\eta_i)$  of the  form 
 \begin{equation}\begin{split}
p_{\e}(\cdot\vert\eta_i)&= 1_{B^{c}}f_{\nu'}(\cdot\vert\eta_i)+\sum_{n=1}^{\infty}1_{B_{n}}(\eta_i)f_{\e,n}(
\cdot\vert\eta_i),
 \quad
{where}\cr
f_{\e,n}(\cdot\vert\eta_i)&=
\left(1-\frac{\e}{n}\right)1_{\{ Y_{\eta_i}\geq 0\} }f_{\nu'}(\cdot\vert\eta_i)+\left(1+\frac{\e}{n}C_
{\eta_i}\right)1_{\{Y_{\eta_i}< 0\} }f_{\nu'}
(\cdot\vert\eta_i).
\end{split}
\end{equation}
Observe that for each $\eta_i\in S'\;$  $p_\e(\cdot|\eta_i)$ is a probability density w.r.t. $\a$ and 
$p_{\e=0}=f_{\nu'}$.

As we did in the proof of part one we introduce a function $\phi:[0,1]\rightarrow \R$ and show that
$\phi$ is decreasing for arguments very close to zero whenever $\nu'(B)>0$.
But this would then imply that $f_{\nu'}$ is not a conditional $\a$-density for a minimizer of $J_{\nu'}$.
Hence for the converse to be true $\nu'(B)=0$, which will then concludes our proof. To formulate this formally
we set $\td\nu^*_\e(d\xi_i)=\nu'(d\eta_i)\a(d\s_i)p_{\e}(\s_i|\eta_i)$ and define
$$\phi(\e):=J_{\nu'}(\td\nu_\e^*)=S\left(\td\nu_\e^*|\a\otimes\nu'\right)+\P(\td\nu_\e^*)-\td\nu_\e^*[\log k].$$
 Now note that the relative entropy takes the form
\begin{equation}
S\left(\td\nu_\e^*|\a\otimes\nu'\right)=\int_{B^{c}}\nu'(d\eta_i)S\left(\td\nu^*(\cdot|\eta_i)\big|\a\right)+
\sum_{n=1}^{\infty} \int_{B_n}\nu'(d\eta_i)S\left(\td\nu^*_\e(\cdot|\eta_i)\big|\a\right).
\end{equation} 
Note that this is a bounded quantity for each $\e\in[0,1]$, since in particular 
$$S\left(\td\nu_\e^*|\a\otimes\nu'\right)\leq 4S\left(\td\nu^*|\a\otimes\nu'\right)+2\log2.$$ Then by the 
Lebesgue's
bounded convergence theorem and the properties of the interaction $\P$ imply that $\phi$ is differentiable 
in the open interval $(0,1)$. 
However, by the convexity of $S$ and once gain the properties of $\P$ one can deduce that 
$$\phi'_+(0)=\lim_{\e\downarrow 0}\frac{\phi(\e)-\phi(0)}{\e}=\frac{d}{d\e}J_{\nu'}(\td\nu_\e^*)\big|_{\e=0}.$$  
Thus the rest of the proof then becomes showing that  $\phi'_+(0)$ has  negative sign whenever $\nu'(B)>0$. Let us now evaluate 
this quantity. Note that because of the form of the perturbed conditional $\a$-density  the part of 
$S\left(\td\nu_\e^*|\a\otimes\nu'\right)$ that will play a
role in determining $\phi'_+(0)$ is the part involving the $B_n$'s. For each $n$ we set
\begin{eqnarray}
S_{n}(\varepsilon)&=&\int_{B_n}\nu'(d\eta_i)S\left(\td\nu^*_\e(\cdot|\eta_i)\big|\a\right)
\end{eqnarray}
 and evaluate 
\begin{equation}\begin{split}
\frac{d}{d\e}S_{n}(\e)\Big\vert_{\e=0}&=\frac{1}{n}\bigg\{\int_{B_{n}}\nu'(d\eta_i)C_{\eta_i} \int_
{Y_{\eta_i}<0}\a(d\s_i) f_{\nu'}
(\s_i\vert\eta_i)\log f_{\nu'}(\s_i\vert\eta_i)\cr
& -\int_{B_{n}} \nu'(d\eta_i)\int_{Y_{\eta_i}\geq0}\a(d\s_i) f_{\nu'}
(\s_i\vert\eta_i)\log f_{\nu'}(\s_i\vert\eta_i)\bigg\}.
\end{split}
\end{equation} 
Next we also evaluate
\begin{equation}\begin{split}
&\frac{d}{d\e}\big(\P(\td\nu^*_\e)-\td\nu^*_\e\left[\log k\right]\big)\Big\vert_{\e=0}\cr
&=\sum_{n=1}^
\infty\frac{1}{n}\bigg\{-\int_{B_{n}}\nu'(d\eta_i) \int_{Y_{\eta_i}\geq0}\a(d\s_i)f_{\nu'}(\s_i\vert\eta_i)
\big(\P^{(1)}
(\pi_1\td\nu^*,\d_{\s_i})
-\log k(\s_i,\eta_i)\big)\cr
&+\int_{B_{n}} \nu'(d\eta_i)C_{\eta_i}\int_{Y_{\eta}<0}\a(d\s_i) f_{\nu'}(\s_i\vert\eta_i)\big
(\P^{(1)}(\pi_1\td\nu^*,\d_{\s_i})
-\log k(\s_i,\eta_i)\big) \bigg\}.
\end{split}
\end{equation} 
It follows from the above considerations and that 
\begin{equation}\begin{split}\label{108}
&\frac{d}{d\e}\big(S(\td\nu^*_\e|\a\otimes\nu')+\P(\td\nu^*_\e)-\td\nu^*_\e\left[\log k\right]\big)
\Big\vert_{\e=0}\cr
&=\sum_{n=1}^
\infty\frac{1}{n}\bigg\{-\int_{B_{n}}\nu'(d\eta_i) \int_{Y_{\eta_i}\geq0}\a(d\s_i)f_{\nu'}
(\s_i\vert\eta_i)\Big( \log f_{\nu'}(\s_i\vert\eta_i)+ \P^{(1)}(\pi_1\td\nu^*,\d_{\s_i})
-\log k(\s_i,\eta_i)\Big)\cr
&+\int_{B_{n}} \nu'(d\eta_i)C_{\eta_i}\int_{Y_{\eta_i}<0}\a(d\s_i)f_{\nu'}(\s_i\vert\eta_i) \Big(\log 
f_{\nu'}(\s_i\vert\eta_i)
+\P^{(1)}(\pi_1\td\nu^*,\d_{\s_i})
-\log k(\s_i,\eta_i)\Big) \bigg\}.
\end{split}
\end{equation}
Furthermore, set  
$$r(\eta_i)=\log \int\a(d\s_i)\exp\left(-\P^{(1)}(\pi_1\td\nu^*,\d_{\s_i})\right)
k(\s_i,\eta_i).$$
Finally, it follows from \eqref{yeta} and \eqref{ceta} after adding and subtracting  $r(\eta_i)$ from the 
integrands of the $\a$-integrals in \eqref{108} that
\begin{equation}\begin{split}
&\frac{d}{d\e}\big(S(\td\nu^*_\e|\a\otimes\nu')+\P(\td\nu^*_\e)-\td\nu^*_\e\left[\log k\right]\big)
\Big\vert_{\e=0}\cr
&=\sum_{n=1}^ 
\infty\frac{1}{n}\bigg\{\int_{B_{n}}\nu'(d\eta_i)r(\eta_i) \int_{Y_{\eta_i}\geq0}\a(d\s_i)f_{\nu'}(\s_i\vert\eta_i)
-\int_{B_{n}}\nu'(d\eta_i) \int_{Y_{\eta_i}\geq0}\a(d\s_i)f_{\nu'}(\s_i\vert\eta_i)Y_{\eta_i}\cr
&-\int_{B_{n}} \nu'(d\eta_i)C_{\eta_i}r(\eta_i)\int_{Y_{\eta_i}<0}\a(d\s_i)f_{\nu'}(\s_i\vert\eta_i)
+\int_{B_{n}} \nu'(d\eta_i)C_{\eta_i}\int_{Y_{\eta_i}<0}\a(d\s_i)f_{\nu'}(\s_i\vert\eta_i)Y_{\eta_i}
\bigg\}\cr
&=\sum_{n=1}^ 
\infty\frac{1}{n}\bigg\{
-\int_{B_{n}}\nu'(d\eta_i) \int_{Y_{\eta_i}\geq0}\a(d\s_i)f_{\nu'}(\s_i\vert\eta_i)Y_{\eta_i}
+\int_{B_{n}} \nu'(d\eta_i)C_{\eta_i}\int_{Y_{\eta_i}<0}\a(d\s_i)f_{\nu'}(\s_i\vert\eta_i)Y_{\eta_i}
\bigg\}\cr
&<0.
\end{split}
\end{equation}
Therefore we have succeeded in showing that 
$$\phi'_+(0)=\lim_{\e\downarrow 0}\frac{\phi(\e)-\phi(0)}{\e}=\lim_{\e\downarrow 0}\frac{J_{\nu'}(\td\nu^*_\e)-
J_{\nu'}(\td\nu^*)}{\e}<0$$ 
whenever $\nu'(B)>0$, contradicting the initial claim that $\td\nu^*$ is a minimizer of $J_{\nu'}$. Hence for 
$\td\nu^*$ to be a minimizer $\nu'(B)=0$. This concludes proof of the claim that the conditional $\a$-densities
$f_{\nu'}$ of the minimizers of $J_{\nu'}$ satisfy the constrained mean-field equation \eqref{cdf} $\a\otimes \nu'$-a.s..
\begin{flushright}
$\Cox$
\end{flushright}


\section*{Acknowledgements }
The authors thank Aernout van Enter, Roberto Fern\'andez, Frank Redig and Wioletta Ruszel  for interesting discussions.

\end{document}